\documentclass[11pt, a4paper]{article}
\usepackage[letterpaper, margin=1.5in]{geometry}
\usepackage{amsmath}
\usepackage{wasysym}
\usepackage{latexsym}
\usepackage{amsfonts}
\usepackage{mathrsfs}
\usepackage{amssymb}
\usepackage{ifsym}
\usepackage{dsfont}
\usepackage[all]{xy}
\usepackage{authblk}

\makeatletter
\def\blfootnote{\xdef\@thefnmark{}\@footnotetext}
\makeatother

\newtheorem{Theorems1}{Theorem}[section]

\newtheorem{Coroll1}[Theorems1]{Corollary}
\newtheorem{Lemma1}[Theorems1]{Lemma}

\newtheorem{Definitions1}[Theorems1]{Definition}
\newtheorem{Examp1}[Theorems1]{Example}

\newenvironment{proof}[1][Proof]{\begin{trivlist}
\item[\hskip \labelsep {\bfseries #1}]}{\end{trivlist}}

\newcommand{\qed}{\nobreak \ifvmode \relax \else
      \ifdim\lastskip<1.5em \hskip-\lastskip
      \hskip1.5em plus0em minus0.5em \fi \nobreak
      \vrule height0.75em width0.5em depth0.25em\fi}

\begin{document}
\title{Initial self-embeddings of models of set theory}
\author[1]{Ali Enayat}
\affil[1]{University of Gothenburg, Gothenburg, Sweden\newline
\texttt{ali.enayat@gu.se}}
\author[2]{Zachiri McKenzie}
\affil[2]{Department of Philosophy, Zhejiang University, Hangzhou, P. R. China\newline
\texttt{zach.mckenzie@gmail.com}}
\maketitle

\begin{abstract}
By a classical theorem of Harvey Friedman (1973), every countable nonstandard model $\mathcal{M}$ of a sufficiently strong fragment of ZF has a proper rank-initial self-embedding $j$, i.e., $j$ is a self-embedding of $\mathcal{M}$  such that $j[\mathcal{M}]\subsetneq\mathcal{M}$, and the ordinal rank of each member of $j[\mathcal{M}]$ is less than the ordinal rank of each element of $\mathcal{M}\setminus j[\mathcal{M}]$.  Here we investigate the larger family of proper \textit{initial-embeddings} $j$ of models $\mathcal{M}$ of fragments of set theory, where the image of $j$ is a transitive submodel of $\mathcal{M}$.  Our results include the following three theorems. In what follows,  $\mathrm{ZF}^-$ is $\mathrm{ZF}$ without the power set axiom; $\mathrm{WO}$ is the axiom stating that every set can be well-ordered; $\mathrm{WF}(\mathcal{M})$ is the well-founded part of $\mathcal{M}$; and  $\Pi^1_\infty\mathrm{-DC}_\alpha$ is the full scheme of dependent choice of length $\alpha$.
\medskip


\noindent \textbf{Theorem A.} \textit{There is an} $\omega$-\textit{standard countable nonstandard model} $\mathcal{M}$ \textit{of} $\mathrm{ZF}^-+\mathrm{WO}$  \textit{that carries no initial self-embedding} $j:\mathcal{M} \longrightarrow \mathcal{M}$ \textit{other than the identity embedding.}

\medskip

\noindent \textbf{Theorem B.} \textit{Every countable $\omega$-nonstandard model $\mathcal{M}$ of} $\mathrm{ZF}$ \textit{is isomorphic to a transitive submodel of the hereditarily countable sets of its own constructible universe} $L^{\mathcal{M}}$.

\medskip
\noindent \textbf{Theorem C.} \textit{The following three conditions are equivalent for a countable nonstandard model} $\mathcal{M}$ \textit{of} $\mathrm{ZF}^{-}+\mathrm{WO}+\forall \alpha\ \Pi^1_\infty\mathrm{-DC}_\alpha$.
\begin{itemize}
\item[(I)] \textit{There is a cardinal in} $\mathcal{M}$ \textit{that is a strict upper bound for the cardinality of each member of} $\mathrm{WF}(\mathcal{M})$.
\item[(II)] $\mathrm{WF}(\mathcal{M})$ \textit{satisfies the powerset axiom}.
\item[(III)] \textit{For all} $n \in \omega$ \textit{and for all} $b \in M$, \textit{there exists a proper initial self-embedding} $j: \mathcal{M} \longrightarrow \mathcal{M}$ \textit{such that} $b \in \mathrm{rng}(j)$ \textit{and} $j[\mathcal{M}] \prec_n \mathcal{M}$.

\end{itemize}


\end{abstract}



\pagebreak
\tableofcontents

\section[Introduction]{Introduction}By a classical theorem of Friedman \cite{fri73}, every countable nonstandard model $\mathcal{M}$ of ZF admits a \textit{proper rank-initial self-embedding}, i.e., an embedding $j: \mathcal{M} \longrightarrow \mathcal{M}$ such that $j[\mathcal{M}]\subsetneq \mathcal{M}$ and the ordinal rank of each member of $\mathcal{M\setminus }j[\mathcal{M}]$ (as computed in $\mathcal{M}$) exceeds the ordinal rank of each member of $j[\mathcal{M}]$ (some authors refer to this situation by saying that $\mathcal{M}$ is a \textit{top} extension of $j[\mathcal{M}]$).  Friedman's work on rank-initial self-embeddings was refined by Ressayre \cite{res87}, who constructed proper rank-initial self-embeddings of models of set theory that pointwise fix any prescribed rank-initial segment of the ambient model determined by an ordinal of the model; and more recently by Gorbow \cite{gor18}, who vastly extended Ressayre's work by carrying out a systematic study of the structure of \textit{fixed point sets} of rank initial self-embeddings of models of set theory.

In another direction, Hamkins \cite{ham13} investigated the family of embeddings $j: \mathcal{M} \longrightarrow \mathcal{N}$, where $\mathcal{M}$ and $\mathcal{N}$ are models of set theory, for which $j[\mathcal{M}]$ is merely required to be a \textit{submodel} of $\mathcal{N}$. The main result of Hamkins' paper shows that, surprisingly, every countable model $\mathcal{M}$ of a sufficiently strong fragment of ZF is embeddable as a submodel of their own constructible universe $L^\mathcal{M}$.

Here we investigate a family of self-embeddings that is wider than the family of rank-initial embeddings, but narrower than the family considered by Hamkins. More specifically, we study \textit{initial} self-embeddings, i.e., embeddings $j:\mathcal{M}\longrightarrow \mathcal{M}$ such that no member of $j[\mathcal{M}]$ gains a new member in the passage from $j[\mathcal{M}]$ to $\mathcal{M}$ (some authors refer to this situation by saying that $\mathcal{M}$ is an \textit{end} extension of $j[\mathcal{M}]$, or that $j[\mathcal{M}]$ is a \textit{transitive} submodel of $\mathcal{M}$). Theorems A, B, and C of the abstract represent the highlights of our results. Theorem A is presented as Theorem \ref{Th:ModelOfZFCminusWithoutSelfEmbedding}; it shows that in contrast with Friedman's aforementioned self-embedding theorem, the theory $\mathrm{ZF^-}$ has countable nonstandard models with no proper initial self-embeddings. Theorem B is presented as Theorem  \ref{Th:MainSelfEmbeddingResultOmegaNonStandard2}; it demonstrates that for $\omega$-nonstandard models, Hamkins' aforementioned theorem can be refined so as to yield a proper \textit{initial} embedding. Finally, Theorem C, which is presented as  Theorem \ref{Th:IsomorphicElementarySubmodelsResult}, gives necessary and sufficient conditions for the existence of proper initial self-embeddings whose images are $\Sigma_n$-elementary in the ambient model; these necessary and sufficient conditions reveal the subtle relationship between the existence of initial self-embeddings of a model $\mathcal{M}$ of set theory and the way in which the well-founded part of $\mathcal{M}$ ``sits" in $\mathcal{M}$.

\section[Preliminaries]{Preliminaries} \label{Sec:Background}

Throughout this paper $\mathcal{L}$ will denote the usual language of set theory whose only nonlogical symbol is the membership relation. Structures will usually be denoted using upper-case calligraphic Roman letters ($\mathcal{M}, \mathcal{N}$, etc.) and the corresponding plain font letter ($M, N$, etc.) will be used to denote the underlying set of that structure.
If $\mathcal{M}$ is an $\mathcal{L}^\prime$-structure where $\mathcal{L}^\prime \supseteq \mathcal{L}$ and $a \in M$, then we will use $a^*$ to denote the set $\{x \in M \mid \mathcal{M} \models (x \in a)\}$ where the background model, $\mathcal{M}$, used in definition of $a^*$ will be clear from the context.

In addition to the L\'{e}vy classes of $\mathcal{L}$-formulae $\Delta_0= \Sigma_1=\Pi_0$, $\Sigma_1$, $\Pi_1$, etc., we will also have cause to consider the Takahashi classes $\Delta_0^\mathcal{P}$, $\Sigma_1^\mathcal{P}$, $\Pi_1^\mathcal{P}$, etc. $\Delta_0^\mathcal{P}$ is the smallest class of $\mathcal{L}$-formulae that contains all atomic formulae, contains all compound formulae formed using the connectives of first-order logic, and is closed under quantification in the form $\mathcal{Q} x \in y$ and $\mathcal{Q} x \subseteq y$ where $x$ and $y$ are distinct variables, and $\mathcal{Q}$ is $\exists$ or $\forall$. The classes $\Sigma_1^\mathcal{P}, \Pi_1^\mathcal{P}$, etc. are defined inductively from the class $\Delta_0^\mathcal{P}$ in the same way that the classes $\Sigma_1, \Pi_1$, etc. are defined from $\Delta_0$. If $\Gamma$ is a collection of formulae and $T$ is a theory, then we will write $\Gamma^T$ for the collection of formulae that are $T$-provably equivalent to a formula in $\Gamma$. If $T$ is an $\mathcal{L}$-theory, then $\Delta_n^T$ is the collection of all $\mathcal{L}$-formulae that are $T$-provably equivalent to both a $\Sigma_n$-formula and a $\Pi_n$-formula. Similarly, $(\Delta_n^\mathcal{P})^T$ is the collection of all $\mathcal{L}$-formulae that are $T$-provably equivalent to both a $\Sigma_n^\mathcal{P}$-formula and a $\Pi_n^\mathcal{P}$-formula.

Let $\mathcal{M}$ and $\mathcal{N}$ be $\mathcal{L}$-structures. We write $\mathcal{M} \equiv \mathcal{N}$ to indicate that $\mathcal{M}$ and $\mathcal{N}$ satisfy the same $\mathcal{L}$-sentences; and write $\mathcal{M} \subseteq \mathcal{N}$ to indicate that $\mathcal{M}$ is a substructure (also referred to as a submodel) of $\mathcal{N}$. If $\Gamma$ is a class of $\mathcal{L}$-formulae, then we will write $\mathcal{M} \prec_\Gamma \mathcal{N}$ if $\mathcal{M} \subseteq \mathcal{N}$ and for every finite tuple $\vec{a} \in M$, $\vec{a}$ satisfies the same $\Gamma$-formulae in both $\mathcal{M}$ and $\mathcal{N}$. In the case that $\Gamma$ is $\Pi_\infty$ (i.e., all $\mathcal{L}$-formulae) or $\Gamma$ is $\Sigma_n$, we will abbreviate this notation by writing $\mathcal{M} \prec \mathcal{N}$ and $\mathcal{M} \prec_n \mathcal{N}$ respectively. If $\mathcal{M} \subseteq \mathcal{N}$ and for all $x \in M$ and $y \in N$,
$$\textrm{if } \mathcal{N} \models (y \in x) \textrm{ then } y \in M,$$
then we say that $\mathcal{N}$ is an \emph{end-extension} of $\mathcal{M}$ (equivalently: $\mathcal{M}$ is an \emph{initial submodel} of $\mathcal{N}$, or  $\mathcal{M}$ is a \emph{transitive submodel} of $\mathcal{N}$) and write $\mathcal{M} \subseteq_e \mathcal{N}$. It is well-known that if $\mathcal{M} \subseteq_e \mathcal{N}$, then $\mathcal{M} \prec_0 \mathcal{N}$. The following is a slight generalisation of the notion of a powerset preserving end-extension that was first studied by Forster and Kaye in \cite{fk91}.

\begin{Definitions1}
Let $\mathcal{M}$ and $\mathcal{N}$ be $\mathcal{L}$-structures. We say that $\mathcal{N}$ is a \textbf{powerset preserving end-extension} of $\mathcal{M}$, and write $\mathcal{M} \subseteq_e^\mathcal{P} \mathcal{N}$ if
\begin{itemize}
\item[(i)] $\mathcal{M} \subseteq_e \mathcal{N}$, and
\item[(ii)] for all $x \in N$ and for all $y \in M$, if $\mathcal{N} \models (x \subseteq y)$, then $x \in M$.
\end{itemize}
\end{Definitions1}

Just as end-extensions preserve $\Delta_0$-properties, powerset preserving end-extensions preserve $\Delta_0^\mathcal{P}$-properties. The following is a slight modification of a result proved in \cite{fk91}:

\begin{Lemma1}
Let $\mathcal{M}$ and $\mathcal{N}$ be $\mathcal{L}$-structures that satisfy Extensionality. If $\mathcal{M} \subseteq_e^\mathcal{P} \mathcal{N}$, then $\mathcal{M} \prec_{\Delta_0^\mathcal{P}} \mathcal{N}$. \Square
\end{Lemma1}

\begin{Definitions1}
Let $\mathcal{M}$ and $\mathcal{N}$ be $\mathcal{L}$-structures. We say that $\mathcal{N}$ is a \textbf{topless powerset preserving end-extension} of $\mathcal{M}$, and write $\mathcal{M} \subseteq_{\mathrm{topless}}^\mathcal{P} \mathcal{N}$ if
\begin{itemize}
\item[(i)] $\mathcal{M} \subsetneq_e^\mathcal{P} \mathcal{N}$, and
\item[(ii)] if $c \in N$ and $c^* \subseteq M$, then $c \in M$.
\end{itemize}
\end{Definitions1}

Let $\Gamma$ be a class of $\mathcal{L}$-formulae. The following define the restriction of some commonly encountered axiom and theorem schemes of $\mathrm{ZFC}$ to formulae in the class $\Gamma$:
\begin{itemize}
\item[]($\Gamma$-Separation) For all $\phi(x, \vec{z}) \in \Gamma$,
$$\forall \vec{z}  \forall w \exists y \forall x(x \in y \iff (x \in w) \land \phi(x, \vec{z})).$$
\item[]($\Gamma$-Collection) For all $\phi(x, y, \vec{z}) \in \Gamma$,
$$\forall \vec{z} \forall w((\forall x \in w) \exists y \phi(x, y, \vec{z}) \Rightarrow \exists c (\forall x \in w)(\exists y \in c) \phi(x, y, \vec{z})).$$
\item[](Strong $\Gamma$-Collection) For all $\phi(x, y, \vec{z}) \in \Gamma$,
$$\forall \vec{z}  \forall w \exists C (\forall x \in w)(\exists y \phi(x, y, \vec{z}) \Rightarrow (\exists y \in C) \phi(x, y, \vec{z})).$$
\item[]($\Gamma$-Foundation) For all $\phi(x, \vec{z}) \in \Gamma$,
$$\forall \vec{z}(\exists x \phi(x, \vec{z}) \Rightarrow \exists y(\phi(y, \vec{z}) \land (\forall x \in y) \neg \phi(x, \vec{z}))).$$
If $\Gamma= \{x \in z\}$ then we will refer to $\Gamma$-Foundation as \emph{Set Foundation}.
\end{itemize}
We will use $\bigcup x \subseteq x$ to abbreviate the $\Delta_0$-formula that says ``$x$ is transitive", i.e., $(\forall y \in x)(\forall z \in y)(z \in x)$. We will also make reference to the axiom of transitive containment ($\mathrm{TCo}$), Zermelo's well-ordering principle ($\mathrm{WO}$), Axiom H and for all $n \in \omega$, the axiom scheme of $\Delta_n$-Separation:
\begin{itemize}
\item[]($\mathrm{TCo}$)
$$\forall x \exists y \left(\bigcup y \subseteq y \land x \subseteq y \right).$$
\item[]($\mathrm{WO}$)
$$\forall x \exists r (r \textrm{ is a well-ordering of } x).$$
\item[](Axiom H)
$$\forall u \exists t \left(\bigcup t \subseteq t \land \forall z(\bigcup z \subseteq z \land |z| \leq |u| \Rightarrow z \subseteq t) \right).$$
\item[]($\Delta_n$-separation) For all $\Sigma_n$-formulae $\phi(x, \vec{z})$ and $\psi(x, \vec{z})$,
$$\forall \vec{z} \ \forall w(\forall x(\phi(x, \vec{z}) \iff \neg \psi(x, \vec{z})) \Rightarrow \exists y \forall x(x \in y \iff (x \in w) \land \phi(x, \vec{z}))).$$
\end{itemize}
For $\alpha$ an ordinal, the {\it $\alpha$-dependent choice scheme} ($\Pi^1_\infty\mathrm{-DC}_\alpha$) is the natural class version of L\'{e}vy's axiom $\mathrm{DC}_\alpha$ \cite{lev64} that generalises Tarski's Dependent Choice Principle by facilitating $\alpha$-sequences of dependent choices.
\begin{itemize}
\item[] ($\Pi^1_\infty\mathrm{-DC}_\alpha$) For all $\mathcal{L}$-formulae $\phi(x, y, \vec{z})$,
$$\forall \vec{z} \left(\begin{array}{c}
\forall g(\forall \gamma \in \alpha)((g \textrm{ is a function})\land (\mathrm{dom}(g)= \gamma) \Rightarrow \exists y \phi(g, y, \vec{z})) \Rightarrow\\
 \exists f (
 (f \textrm{ is a function})\land (\mathrm{dom}(f)= \alpha) \land
 (\forall \beta \in \alpha)\phi(f \upharpoonright \beta, f(\beta), \vec{z}))
\end{array}\right).$$
\end{itemize}

We will have cause to consider the following subsystems of $\mathrm{ZFC}$:
\begin{itemize}
\item $\mathbf{M}^-$ is the $\mathcal{L}$-theory with axioms: Extensionality, Emptyset, Pair, Union, Infinity, $\mathrm{TCo}$, $\Delta_0$-Separation and Set Foundation.
\item $\mathbf{M}$ is obtained from $\mathbf{M}^-$ by adding the powerset axiom.
\item $\mathrm{Mac}$ is obtained from $\mathbf{M}$ by adding the axiom of choice.
\item $\mathrm{KPI}$ is obtained from $\mathbf{M}^-$ by adding $\Delta_0$-Collection and $\Pi_1$-Foundation.
\item $\mathrm{KP}$ is obtained from $\mathrm{KPI}$ by removing the axiom of infinity.
\item $\mathrm{KP}^\mathcal{P}$ is obtained from $\mathbf{M}$ by adding $\Delta_0^\mathcal{P}$-Collection and $\Pi_1^\mathcal{P}$-Foundation.
\item $\mathrm{MOST}$ is obtained from $\mathrm{Mac}$ by adding $\Sigma_1$-Separation and $\Delta_0$-Collection.
\item $\mathrm{ZF}^-$ is obtained by adding $\Pi_\infty$-Collection to $\mathrm{KPI}$.



\end{itemize}
The theories $\mathbf{M}$, $\mathrm{KPI}$, $\mathrm{KP}$ and $\mathrm{KP}^\mathcal{P}$ are studied in \cite{mat01}. In contrast with the version of Kripke-Platek Set Theory studied in \cite{fri73, bar75}, which includes $\Pi_\infty$-Foundation, we follow \cite{mat01}, by only including $\Pi_1$-Foundation in the theories $\mathrm{KP}$ and $\mathrm{KPI}$, and only including $\Pi_1^\mathcal{P}$-Foundation in the theory $\mathrm{KP}^\mathcal{P}$. The theory $\mathrm{KPI}$, as defined here, plays a key role in \cite{flw16}, where it is referred to as $\mathrm{KP}^- + \mathrm{infinity} + \Pi_1$-Foundation.

The results of \cite{zar96} and, more recently, \cite{ght16} highlight the importance of axiomatising $\mathrm{ZF}^-+\mathrm{WO}$ using the collection scheme ($\Pi_\infty$-Collection) instead of the replacement scheme. The strength of Zermelo's well-order principle $\mathrm{WO}$ in the $\mathrm{ZF}^-$ context is revealed in \cite{zar82}, which shows that, in the absence of the powerset axiom, the statement that {\it every set of nonempty sets has a choice function} does not imply $\mathrm{WO}$.\footnote{Zarach credits Z. Szczepaniak with first finding a model of $\mathrm{ZF}^-$ in which every set of nonempty sets has a choice function, but in which $\mathrm{WO}$ fails.}

Let $\mathrm{ZF}^-+ \mathrm{GWO}(R)$ be the extension of $\mathrm{ZF}^-+\mathrm{WO}$ obtained as follows: introduce a new binary relation symbol, $R$, to the language of set theory, and then add an axiom asserting that $R$ is a bijection between the universe and the class of ordinals (a global well-order), and also extend the schemes of separation and collection so as to ensure that formulae mentioning $R$ can be used. As a consequence of a result of Flanagan \cite[Theorem 7.1]{fla75}, $\mathrm{ZF}^-+ \mathrm{GWO}(R)$ is a conservative extension of the theory $\mathrm{ZF}^-+\mathrm{WO}+\forall \alpha \ \Pi^1_\infty\mathrm{-DC}_\alpha$.  Recent work of S. Friedman, Gitman and Kanovei \cite{fgk19} shows that $\Pi^1_\infty\mathrm{-DC}_{\omega}$ is independent of $\mathrm{ZF}^-+\mathrm{WO}$.

Next we record the following useful relationships between fragments of Collection, Separation and Foundation over the base theory $\mathbf{M}^-$.


\begin{Lemma1} \label{basicimplications}
Let $\Gamma$ be a class of $\mathcal{L}$-formulae, and $n \in \omega$.
\begin{enumerate}
\item In the presence of $\mathbf{M}^-$,  $\Pi_n\textrm{-Separation}$ is equivalent to $\Sigma_n\textrm{-Separation}.$
\item $\mathbf{M}^-+\Gamma\textrm{-Separation} \vdash \Gamma\textrm{-Foundation}.$

\item $\mathbf{M}^{-}+\Pi_n\textrm{-Collection} \vdash \Sigma_{n+1}\textrm{-Collection}.$


\item \cite[Lemma 4.13]{flw16} $\mathbf{M}^-+\Pi_n\textrm{-Collection} \vdash \Delta_{n+1}\textrm{-Separation}$.
\item \cite[Lemma 2.5]{mck19} In the presence of $\mathbf{M}^-$, $\Pi_n\textrm{-Collection}+\Sigma_{n+1}\textrm{-Separation}$ is equivalent to $\textrm{Strong }\Pi_n\textrm{-Collection}$.
\end{enumerate}
\end{Lemma1}

As indicated by the following well-known result, over the theory $\mathbf{M}^-$, $\Pi_n$-Collection implies that the classes $\Pi_{n+1}$ and $\Sigma_{n+1}$ are essentially closed under bounded quantification (part (3) of Lemma \ref{basicimplications} is used in the proof).

\begin{Lemma1}
Let $\phi(x, \vec{z})$ be a $\Sigma_{n+1}$-formula and let $\psi(x, \vec{z})$ be a $\Pi_{n+1}$-formula. The theory $\mathbf{M}^-+\Pi_n\textrm{-Collection}$ proves that $(\forall x \in y) \phi(x, \vec{z})$ is equivalent to a $\Sigma_{n+1}$-formula and $(\exists x \in y)\psi(x, \vec{z})$ is equivalent to a $\Pi_{n+1}$-formula.
\end{Lemma1}

\begin{Definitions1}
A transitive set $M$ is said to be \textbf{admissible} if $\langle M, \in\rangle \models \mathrm{KP}$.
\end{Definitions1}

The theory $\mathrm{KPI}$ and its variants that include the scheme of full class foundation have been widely studied \cite{fri73, bar75, mat01, flw16}. One appealing feature of this theory is the fact that it is strong enough to carry out many of the fundamental set-theoretic constructions such as defining set-theoretic rank, proving the existence of transitive closures, defining satisfaction and constructing G\"{o}del's $L$ hierarchy.
\begin{itemize}
\item For all sets $X$, we use $\mathrm{TC}(X)$ to denote the $\subseteq$-least transitive set with $X$ as a subset. The theory $\mathrm{KPI}$ proves that the function $X \mapsto \mathrm{TC}(X)$ is total. Moreover, the proof of \cite[Proposition 1.29]{mat01} shows that the formulae ``$x = \mathrm{TC}(y)$" and ``$x \in \mathrm{TC}(y)$" with free variables $x$ and $y$ are $\Delta_1^\mathrm{KP}$, and ``$x = \mathrm{TC}(y)$" is also $\Delta_0^\mathcal{P}$.
\item The theory $\mathrm{KP}$ is capable of defining and proving the totality of the rank function $\rho$ satisfying
$$\rho(a)= \sup \left\{ \rho(b)+1: b\in a \right\}.$$
The formula ``$\rho(x)=y$" with free variables $x$ and $y$ is $\Delta_1^{\mathrm{KP}}$ \cite[Theorem 1.5]{fri73}.
\item As verified in \cite[Section III.1]{bar75}, satisfaction in set structures is definable in $\mathrm{KPI}$. In particular, if $\mathcal{N}$ is a set structure in a model $\mathcal{M}$ of $\mathrm{KPI}$, $\mathcal{L}(\mathcal{N})$ is the language of $\mathcal{N}$, $\vec{a}$ is an $\mathcal{M}$-finite sequence of members of $\mathcal{N}$, and $\phi$ is an $\mathcal{L}(\mathcal{N})$-formula in the sense of $\mathcal{M}$ whose arity agrees with the length of $\vec{a}$, then ``$\mathcal{M} \models \phi[\vec{v}/\vec{a}]$" is definable in $\mathcal{M}$ by a formula that is $\Delta_1^{\mathrm{KPI}}$.
\item As shown in \cite[Chapter II]{bar75} the theory $\mathrm{KPI}$ is capable of constructing the levels of G\"{o}del's $L$ hierarchy. The following operation can be defined using a formula for satisfaction for set structures in $\mathrm{KPI}$: for all sets $X$,
$$\mathrm{Def}(X)= \{Y \subseteq X \mid Y \textrm{ is a definable subclass of } \langle X, \in \rangle\}.$$
The levels of the $L$ hierarchy are then recursively defined by:
$$L_0= \emptyset \textrm{ and } L_\alpha= \bigcup_{\beta < \alpha} L_\beta \textrm{ if }\alpha \textrm{ is a limit ordinal},$$
$$L_{\alpha+1}= L_\alpha \cup \mathrm{Def}(L_\alpha),$$
$$L= \bigcup_{\alpha \in \mathrm{Ord}} L_\alpha.$$
The function $\alpha \mapsto L_\alpha$ is total and $\Delta_1^{\mathrm{KP}}$. As usual, we will use $V=L$ to abbreviate the axiom that says that every set is a member of some $L_\alpha$, i.e., $\forall x \exists \alpha (x \in L_\alpha)$.
\end{itemize}

The fact that $\mathrm{KPI}$ can express satisfaction in set structures can be used, in this theory, to express satisfaction for $\Delta_0$-formulae in the universe via the definition below.

\begin{Definitions1} \label{Df:Delta0Satisfaction}
The formula $\mathrm{Sat}_{\Delta_0}(q, x)$ is defined as
$$\begin{array}{c}
(q \in \omega) \land (q= \ulcorner \phi(v_1, \ldots, v_m) \urcorner \textrm{ where } \phi \textrm{ is } \Delta_0) \land (x= \langle x_1, \ldots, x_m \rangle) \land\\
\exists N \left( \bigcup N \subseteq N \land (x_1, \ldots, x_m \in N) \land (\langle N, \in \rangle \models \phi[x_1, \ldots, x_m]) \right)
\end{array}.$$
\end{Definitions1}

The absoluteness of $\Delta_0$ properties between transitive structures and the universe, and the availability of $\mathrm{TCo}$ in $\mathrm{KPI}$ imply that the formula $\mathrm{Sat}_{\Delta_0}$ is equivalent, in the theory $\mathrm{KPI}$, to the formula
$$\begin{array}{c}
(q \in \omega) \land (q= \ulcorner \phi(v_1, \ldots, v_m) \urcorner \textrm{ where } \phi \textrm{ is } \Delta_0) \land (x= \langle x_1, \ldots, x_m \rangle) \land\\
\forall N \left( \bigcup N \subseteq N \land (x_1, \ldots, x_m \in N) \Rightarrow (\langle N, \in \rangle \models \phi[x_1, \ldots, x_m]) \right)
\end{array}.$$
Therefore, the fact that ``$\langle N, \in \rangle \models \phi[\cdots]$" is $\Delta_1^{\mathrm{KPI}}$ implies that $\mathrm{Sat}_{\Delta_0}(q, x)$ is also $\Delta_1^{\mathrm{KPI}}$, and $\mathrm{Sat}_{\Delta_0}(q, x)$ expresses satisfaction for $\Delta_0$-formulae in the theory $\mathrm{KPI}$. We can now inductively define formulae $\mathrm{Sat}_{\Sigma_n}(q, x)$ and $\mathrm{Sat}_{\Pi_n}(q, x)$ that express satisfaction for formulae in the classes $\Sigma_n$ and $\Pi_n$.

\begin{Definitions1}
The formulae $\mathrm{Sat}_{\Sigma_n}(q, x)$ and $\mathrm{Sat}_{\Pi_n}(q, x)$ are defined recursively for $n>0$. $\mathrm{Sat}_{\Sigma_{n+1}}(q, x)$ is defined as the  formula
$$\exists \vec{y} \exists k \exists b \left( \begin{array}{c}
(q= \ulcorner\exists \vec{u} \phi(\vec{u}, v_1, \ldots, v_l)\urcorner \textrm{ where } \phi \textrm{ is } \Pi_n)\land (x= \langle x_1, \ldots, x_l \rangle)\\
\land (b= \langle \vec{y}, x_1, \ldots, x_l \rangle) \land (k= \ulcorner \phi(\vec{u}, v_1, \ldots, v_l) \urcorner) \land \mathrm{Sat}_{\Pi_n}(k, b)
\end{array}\right);$$
and  $\mathrm{Sat}_{\Pi_{n+1}}(q, x)$ is defined as the formula
$$\forall \vec{y} \forall k \forall b \left( \begin{array}{c}
(q= \ulcorner\forall \vec{u} \phi(\vec{u}, v_1, \ldots, v_l) \urcorner \textrm{ where } \phi \textrm{ is } \Sigma_n)\land (x= \langle x_1, \ldots, x_l \rangle)\\
\land ((b= \langle \vec{y}, x_1, \ldots, x_l \rangle) \land (k= \ulcorner\phi(\vec{u}, v_1, \ldots, v_l)\urcorner) \Rightarrow \mathrm{Sat}_{\Sigma_n}(k, b))
\end{array}\right).$$
\end{Definitions1}

\begin{Theorems1} \label{Complexityofpartialsat} Suppose $n \in \omega$ and $m=\max \{ 1, n \}$. The formula $\mathrm{Sat}_{\Sigma_n}(q, x)$ (respectively $\mathrm{Sat}_{\Pi_n}(q, x)$) is $\Sigma_n^{\mathrm{KPI}}$ ($\Pi_n^{\mathrm{KPI}}$, respectively). Moreover, $\mathrm{Sat}_{\Sigma_n}(q, x)$ (respectively $\mathrm{Sat}_{\Pi_n}(q, x)$) expresses satisfaction for $\Sigma_n$-formulae ($\Pi_n$-formulae, respectively) in the theory $\mathrm{KPI}$, i.e., if $\mathcal{M} \models\mathrm{KPI}$, $\phi(v_1,\ldots,v_k)$ is a $\Sigma_n$-formula, and $x_1,\ldots,x_k$ are in $M$, then for $q = \ulcorner   \phi( v_1, \ldots, v_k) \urcorner$, $\mathcal{M}$ satisfies the universal generalization of the following formula:

$$  x= \langle x_1, \ldots,x_k \rangle \Rightarrow \left( \phi(x_1,\ldots,x_k) \leftrightarrow    \mathrm{Sat}_{\Sigma_{n}}(q, x) \right).$$

\end{Theorems1}

The following result appears in  \cite[Theorem 3.8]{flw16}.

\begin{Lemma1} \label{Th:SchroderBernsteinInKP}
(Friedman, Li, Wong) The theory $\mathrm{KP}$ proves the Schr\"{o}der-Bernstein Theorem, i.e.,  $\mathrm{KP}$ proves that if $A$ and $B$ are sets such that $|A| \leq |B|$ and $|B| \leq |A|$, then $|A|=|B|$.
\end{Lemma1}

The following theorem highlights the important fact that the $\Sigma_1^\mathcal{P}$-Recursion Theorem is provable in the theory $\mathrm{KP}^\mathcal{P}$ \cite[Theorem 6.26]{mat01}.

\begin{Theorems1} \label{Th:Sigma1PRecursionTheorem}
($\mathrm{KP}^\mathcal{P}$) Let $G$ be a $\Sigma_1^\mathcal{P}$-definable class. If $G$ is a total function, then there exists a $\Sigma_1^\mathcal{P}$-definable total class function $F$ such that for all $x$, $F(x)= G(F\upharpoonright x)$.
\end{Theorems1}

\begin{Definitions1} \label{V_alpha}
We write ``$V_\alpha$ exists'' as an abbreviation for the sentence expressing that $\alpha$ is an ordinal, and there is a function $f$ whose domain is
$\alpha+1$ that satisfies the following conditions (1) through (3) below.
\begin{enumerate}
    \item $f(0)=\emptyset$.
    \item $\forall \beta<\alpha \left((\beta \textrm{ is a limit ordinal}) \Rightarrow f(\beta)= \bigcup_{\xi < \beta} f(\xi) \right)$.
    \item $(\forall \beta \in \mathrm{dom}(f))(\forall y )(y \in f(\beta+1) \iff y \subseteq f(\beta))$.
\end{enumerate}
\end{Definitions1}

Note that under Definition \ref{V_alpha}, if  $V_\alpha$ exists, then $V_\beta$ exists for all $\beta<\alpha$. The following consequence of the $\Sigma_1^\mathcal{P}$-Recursion Theorem is Proposition 6.28 of \cite{mat01}.

\begin{Coroll1} \label{Th:RanksInKPP}
The theory $\mathrm{KP}^\mathcal{P}$ proves that for all ordinals $\alpha$, $V_\alpha$ exists. Note that in particular, this theory proves that for all ordinals $\alpha$, there is a function $f$ with domain $\alpha$ such that for all $\beta \in \alpha$,  $f(\beta) = V_\beta$. \Square
\end{Coroll1}

Section 3 of \cite{mat01} contains the verification of the following lemma.

\begin{Lemma1} \label{Th:MOSTisMacPlusH}
$\mathrm{MOST}$ is the theory $\mathrm{Mac}+\textrm{Axiom }\mathrm{H}$.
\end{Lemma1}

We also record the following consequence of $\mathrm{MOST}$ that are proved in \cite[Section 3]{mat01}:

\begin{Lemma1} \label{Th:ConsequencesOfMOST}
The theory $\mathrm{MOST}$ proves
\begin{itemize}
\item[(i)] every well-ordering is isomorphic to an ordinal,
\item[(ii)] every well-founded extensional relation is isomorphic to a transitive set,
\item[(iii)] for all cardinals $\kappa$, $\kappa^+$ exists, and
\item[(iv)] for all cardinals $\kappa$, the set $H_\kappa= \{x \mid |\mathrm{TC}(x)| < \kappa \}$ exists.
\end{itemize}
\end{Lemma1}

The following result is \cite[Lemma 3.3]{ekm18} combined with the refinement of a theorem due to Takahashi proved in \cite[Proposition Scheme 6.12]{mat01}:


\begin{Lemma1} \label{Th:HCutsSatisfyPi1Collection}
If $\mathcal{M}, \mathcal{N} \models \mathrm{MOST}$ and $\mathcal{N} \subseteq_{\mathrm{topless}}^\mathcal{P} \mathcal{M}$, then $\mathcal{N} \models \Pi_1\textrm{-Collection}$.
\end{Lemma1}


We next recall a remarkable absoluteness phenomenon unveiled by  L\'{e}vy \cite{lev64}, which shows that, provably in ZF, $H^{L}_{\aleph_1}$ (i.e., the collection of sets that are hereditarily countable, as computed in the constructible universe) is a $\Sigma_1$-elementary submodel of the universe of sets.\footnote{The proof of this result relies heavily on  the venerable Shoenfield Absoluteness Theorem. The original proof by L\'{e}vy of Theorem \ref{Levy-Shoenfield} presented this result as a theorem of ZF + DC (where DC here is axiom of dependent choice of length $\omega$).  As pointed out by Kunen, DC can be eliminated by a forcing-and-absoluteness stratagem (see  page 55 of \cite{bar71}). Later Barwise and Fischer gave a direct forcing-free proof in ZF \cite{barwise-fischer70}.}

\begin{Theorems1} \label{Levy-Shoenfield}
(L\'{e}vy-Shoenfield Absoluteness) Let $\theta(x,y)$ be a $\Sigma_1^{\mathrm{ZF}}$-formula with no free variables except $x$ and $y$, then the universal generalization of the following formula is provable in $\mathrm{ZF}$

$$ \left(y \in H^{L}_{\aleph_1} \land \exists x\ \theta(x,y) \right)\Rightarrow \exists x(x\in H^{L}_{\aleph_1} \land \theta(x,y)).$$

\end{Theorems1}

The L\'{e}vy-Shoenfield Absoluteness Theorem readily implies the following corollary that shows that the $\Sigma_1$-theory of every model of ZF coincides with the $\Sigma_1$-theory of $H_{\aleph_1}$ of the constructible universe of the same model.

\begin{Coroll1} \label{cor. of Levy-Shoenfield}  Let $\delta(\vec{x})$ be a $\Delta_0^{\mathrm{ZF}}$-formula, and $\mathcal{M}$ be a model of $\mathrm{ZF}$. Then

$$\mathcal{M} \models \exists \vec{x}\ \delta(\vec{x}) \Longleftrightarrow H^{L^{\mathcal{M}}}_{\aleph_1} \models \exists \vec{x}\ \delta(\vec{x}).$$
\end{Coroll1}

Any model of $\mathrm{KPI}$ comes equipped with its well-founded part that consists of all sets in this structure whose rank is a standard ordinal, as indicated by the following definition.


\begin{Definitions1} \label{Df:StandardPart}
Let $\mathcal{M} \models \mathrm{KP}$. The \textbf{well-founded part} or \textbf{standard part} of $\mathcal{M}$, denoted $\mathrm{WF}(\mathcal{M})$, is the substructure of $\mathcal{M}$ with underlying set
$$\mathrm{WF}(M)= \{x \in M \mid \neg \exists f(f: \omega \longrightarrow M\land(\forall n \in \omega)(f(0)= x \land f(n+1)\in^\mathcal{M} f(n)))\}.$$
If $\mathrm{WF}(\mathcal{M})\neq \mathcal{M}$, then we say that $\mathcal{M}$ is {\bf nonstandard}. The \textbf{standard ordinals} of $\mathcal{M}$, denoted $\mathrm{o}(\mathcal{M})$, is the substructure of $\mathcal{M}$ with underlying set $\mathrm{o}(M)= \mathrm{WF}(M)\cap \mathrm{Ord}^\mathcal{M}$. If $\omega^\mathcal{M} \in \mathrm{o}(M)$, then we say that $\mathcal{M}$ is $\mathbf{\omega}$\textbf{-standard}; otherwise $\mathcal{M}$ is said to be $\mathbf{\omega}$\textbf{-nonstandard}. Mostowski's Collapsing Lemma ensures that both $\mathrm{o}(\mathcal{M})$ and $\mathrm{WF}(\mathcal{M})$ are isomorphic to transitive sets. In particular, $\mathrm{o}(\mathcal{M})$ is isomorphic to an ordinal that is called the \textbf {standard ordinal} of $\mathcal{M}$.
\end{Definitions1}

The following definition generalises the notion of standard system that plays an important role in the study of models of arithmetic.

\begin{Definitions1}
Let $\mathcal{M} \models \mathrm{KPI}$. The \textbf{standard system} of $\mathcal{M}$ is the set
$$\mathrm{SSy}(\mathcal{M})= \{y^* \cap \mathrm{WF}(\mathcal{M}) \mid y \in M\}.$$
If $A \in \mathrm{SSy}(\mathcal{M})$ and $y \in M$ is such that $A= y^* \cap \mathrm{WF}(\mathcal{M})$, then we say that $y$ {\bf codes} $A$.
\end{Definitions1}


\begin{Definitions1}
Let $\mathcal{M}$ and $\mathcal{N}$ be $\mathcal{L}$-structures. An \textbf{embedding} of $\mathcal{M}$ into $\mathcal{N}$ is an injection $j: M \longrightarrow N$ such that for all $x, y \in M$,
$$\mathcal{M} \models x \in y \textrm{ if and only if } \mathcal{N} \models j(x) \in j(y).$$
Note that we will often write $j: \mathcal{M} \longrightarrow \mathcal{N}$ to indicate that $j$ is an embedding of $\mathcal{M}$ into $\mathcal{N}$. If $j: \mathcal{M} \longrightarrow \mathcal{N}$ is an embedding of $\mathcal{M}$ into $\mathcal{N}$, then we write $j[\mathcal{M}]$ for the substructure of $\mathcal{N}$ whose underlying set is $\mathrm{rng}(j)$.
\end{Definitions1}

\begin{Definitions1}
Let $\mathcal{M}$ be an $\mathcal{L}$-structure and let $j: \mathcal{M} \longrightarrow \mathcal{M}$ be an embedding of $\mathcal{M}$ into $\mathcal{M}$. The \textbf{fixed point set} of $j$ is the set $\mathrm{Fix}(j)= \{ x \in M \mid j(x)=x\}$.
\end{Definitions1}

\begin{Definitions1}
Let $\mathcal{M}$ and $\mathcal{N}$ be $\mathcal{L}$-structures. Let $j: \mathcal{M} \longrightarrow \mathcal{N}$ be an embedding of the structure $\mathcal{M}$ into $\mathcal{N}$. We say that $j$ is an \textbf {initial embedding} if $j[\mathcal{M}] \subseteq_e \mathcal{N}$. We say that $j$ is a \textbf {$\mathcal{P}$-initial embedding} if $j[\mathcal{M}] \subseteq_e^\mathcal{P} \mathcal{N}$. If $j: \mathcal{M} \longrightarrow \mathcal{M}$ is a ($\mathcal{P}$-) initial embedding with $j[\mathcal{M}]\neq \mathcal{M}$, then we say that $j$ is a \textbf{proper ($\mathcal{P}$-) initial self-embedding} of $\mathcal{M}$.
\end{Definitions1}

Next, we take advantage of the rank function available in $\mathrm{KP}$ to define the notion of rank extension, and the notion of rank-initial embedding.

\begin{Definitions1}
Let $\mathcal{M}$ and $\mathcal{N}$ be $\mathcal{L}$-structures with $\mathcal{M} \subseteq_e \mathcal{N}$ and $\mathcal{N} \models \mathrm{KP}$. We say that $\mathcal{N}$ is a \textbf{rank extension} of $\mathcal{M}$ if for all $\alpha \in \mathrm{Ord}^\mathcal{M}$, if $x \in N$ and $\mathcal{N} \models \rho(x)=\alpha$, then $x \in M$.
\end{Definitions1}


\begin{Definitions1}
Let $\mathcal{M}$ and $\mathcal{N}$ be $\mathcal{L}$-structures that satisfy $\mathrm{KP}$. Let $j: \mathcal{M} \longrightarrow \mathcal{N}$ be an embedding of the structure $\mathcal{M}$ into $\mathcal{N}$. We say $j$ is a \textbf{rank-initial embedding} if $j$ is an initial embedding and $\mathcal{M}$ is a rank extension of $j[\mathcal{M}]$. If $j: \mathcal{M} \longrightarrow \mathcal{M}$ is a rank-initial embedding with $j[\mathcal{M}]\neq \mathcal{M}$, then we say that $j$ is a \textbf{proper rank-initial self-embedding} of $\mathcal{M}$.
\end{Definitions1}

Note that a rank-initial embedding $j: \mathcal{M} \longrightarrow \mathcal{N}$, where $\mathcal{M}, \mathcal{N} \models \mathrm{KP}$, is also $\mathcal{P}$-initial. The following result of Gorbow \cite[Corollary 4.6.12]{gor18} shows that if the source and target model of a $\mathcal{P}$-initial embedding both satisfy $\mathrm{KP}^\mathcal{P}$, then this embedding is also rank-initial.

\begin{Lemma1} \label{Th:PInitialImpliesRankInitialKPP}
Let $\mathcal{M}$ and $\mathcal{N}$ be models of $\mathrm{KP}^\mathcal{P}$. If $j: \mathcal{M} \longrightarrow \mathcal{N}$ is a $\mathcal{P}$-initial embedding, then $j$ is a rank-initial embedding.
\end{Lemma1}

Note that, in any model of $\mathrm{ZFC}$, $L$ is a powerset preserving end-extension of $H_{\aleph_\omega}^L$ and $H_{\aleph_\omega}^L$ satisfies $\mathrm{MOST}+\Pi_\infty\textrm{-Separation}$. This example shows that the assumption that $\mathcal{M}$ satisfies $\Delta_0^\mathcal{P}\textrm{-Collection}$ in Lemma \ref{Th:PInitialImpliesRankInitialKPP} can not relaxed to $\Delta_0\textrm{-Collection}$ even in the presence of the full scheme of separation.

H. Friedman's seminal \cite{fri73} pioneered the study of rank-initial self-embeddings of $\mathrm{KP}^\mathcal{P}+\Pi_\infty\textrm{-Foundation}$. His work was refined and extended by Ressayre \cite{res87}, and more recently by Gorbow \cite{gor18}. The following theorem of Gorbow guarantees the existence of proper rank-initial self-embeddings of countable nonstandard models of an extension of $\mathrm{KP}^\mathcal{P}$. Gorbow's theorem refines  \cite[Theorem 4.3]{fri73}, and is a consequence of results proved in \cite[Section 5.2]{gor18}.


\begin{Theorems1} (Gorbow) \label{Th_Gorbow}
Every countable nonstandard model $\mathcal{M}$ of $\mathrm{KP}^\mathcal{P}+\Sigma_1^\mathcal{P}\textrm{-Separation}$ has a proper rank-initial self-embedding. Moreover, given any $\alpha \in \mathrm{Ord}^\mathcal{M}$ there exists a proper rank-initial self-embedding $j$ of $ \mathcal{M}$ that fixes every element of $(\mathrm{V}_{\alpha}^\mathcal{M})^*$.

\end{Theorems1}

We also note the following self-embedding theorem that is readily obtained by putting \cite[Theorem 5.6]{ekm18} together \cite[Proposition Scheme 6.12]{mat01}.


\begin{Theorems1} \label{Th:PInitialSelfEmbeddingOfRecursivelySaturatedModels}
Every countable recursively saturated model of $\mathrm{MOST}+\Pi_1\textrm{-Collection}$ has a proper $\mathcal{P}$-initial self-embedding.
\end{Theorems1}


\section[The well-founded part]{The well-founded part}

In this section we present results about well-founded parts of models of set theory that are relevant to the proofs the main result of this paper. H. Friedman's \cite{fri73} systematically studied the structure of the well-founded part of a nonstandard model of Kripke-Platek Set Theory and \cite[Theorem 2.1]{fri73} showed that such a well-founded part must be isomorphic to an admissible set. This result is also a consequence of \cite[Lemma 8.4]{bar75}. As we mentioned earlier, the versions of Kripke-Platek Set Theory studied in \cite{bar75} and \cite{fri73} include $\Pi_\infty$-Foundation. An examination of these proofs reveals that the well-founded part of a model of $\mathrm{KPI}+\Sigma_1\textrm{-Foundation}$ is isomorphic to an admissible set. Before proving this, we will first verify in the lemma below that any nonstandard model of $\mathrm{KPI}$ is a topless powerset preserving end-extension of its well-founded part.

\begin{Lemma1} \label{Th:WellfoundedPartTopless}
Let $\mathcal{M}$ be a nonstandard model of $\mathrm{KPI}$. Then
$$\mathrm{WF}(\mathcal{M}) \subseteq_{\mathrm{topless}}^\mathcal{P} \mathcal{M}.$$
\end{Lemma1}

\begin{proof}
It follows immediately from Definition \ref{Df:StandardPart} that $\mathrm{WF}(\mathcal{M}) \subseteq_{e}^\mathcal{P} \mathcal{M}$. The fact that $\mathcal{M}$ is nonstandard immediately means that $M \neq \mathrm{WF}(M)$. Let $c \in M$ with $c^* \subseteq \mathrm{WF}(M)$. Suppose, for a contradiction, that $f: \omega \longrightarrow M$ witness the fact that $c \notin \mathrm{WF}(M)$. But, $f(0)= c$ and $\mathcal{M} \models (f(1) \in c)$, so $f(1) \in \mathrm{WF}(M)$. Define $g: \omega \longrightarrow M$ by: for all $n \in \omega$ with $n \geq 1$, $g(n-1)=f(n)$. Now, $g$ witness the fact that $f(1) \notin \mathrm{WF}(M)$, which is a contradiction. This shows that $\mathrm{WF}(\mathcal{M}) \subseteq_{\mathrm{topless}}^\mathcal{P} \mathcal{M}$.
\Square
\end{proof}

\begin{Lemma1} \label{Th:StandardPartClosedUnderRank}
Let $\mathcal{M}$ be a nonstandard model of $\mathrm{KPI}$. Then $\mathrm{WF}(\mathcal{M})$ satisfies Extensionality, Emptyset, Pair, Union, $\Delta_0$-Separation and $\Pi_\infty$-Foundation. Moreover, for all $x \in M$,
$$x \in \mathrm{WF}(M) \textrm{ if and only if } \rho^\mathcal{M}(x) \in \mathrm{WF}(M).$$
\end{Lemma1}

\begin{proof}
The fact $\mathrm{WF}(\mathcal{M}) \subseteq_e^\mathcal{P} \mathcal{M}$ implies that $\mathrm{WF}(\mathcal{M})$ satisfies Extensionality, Emptyset, Pair, Union, and $\Delta_0$-Separation. The fact that $\mathrm{WF}(\mathcal{M})$ is well-founded ensures that $\Pi_\infty$-Foundation holds in $\mathrm{WF}(\mathcal{M})$. Now, to prove the last statement, consider
$$a= \{x \in \mathrm{WF}(M) \mid \rho^\mathcal{M}(x) \notin \mathrm{WF}(M) \}.$$
Suppose that $a \neq \emptyset$ and let $z \in a$ be a $\in^\mathcal{M}$-least member of $a$. Working inside $\mathcal{M}$, consider
$$b= \{ \rho(y)+1 \mid y \in z\}.$$
Next we observe that $b$ is a set and $b^* \subseteq \mathrm{WF}(M)$,  so $b \in \mathrm{WF}(M)$. It now follows that $\rho^\mathcal{M}(z)= (\sup b)^\mathcal{M} \in \mathrm{WF}(M)$, which is a contradiction. This shows that if $x \in \mathrm{WF}(M)$, then $\rho^\mathcal{M}(x) \in \mathrm{WF}(M)$. Conversely, consider
$$c= \{\rho^\mathcal{M}(x) \mid x \notin \mathrm{WF}(M)\}.$$
Suppose, for a contradiction, that $c \cap \mathrm{WF}(M)\neq \emptyset$. Let $\alpha \in c$ be least and let $z \in M$ with $z \notin \mathrm{WF}(M)$ be such that $\mathcal{M} \models \rho(z)= \alpha$. Since $\alpha$ is least, $z^* \subseteq \mathrm{WF}(M)$. Since $\mathrm{WF}(\mathcal{M}) \subseteq_{\mathrm{topless}}^\mathcal{P} \mathcal{M}$, this implies that $z \in \mathrm{WF}(M)$, which is a contradiction. \Square
\end{proof}

We next verify that under the additional assumption that the nonstandard model of $\mathrm{KPI}$ satisfies $\Sigma_1$-Foundation, the well-founded part also satisfies $\Delta_0$-Collection.

\begin{Theorems1} \label{Th:WellFoundedPartIsAdmissible}
Let $\mathcal{M}$ be a nonstandard model of $\mathrm{KPI}+\Sigma_1\textrm{-Foundation}$. Then $\mathrm{WF}(\mathcal{M})$ is isomorphic to an admissible set.
\end{Theorems1}

\begin{proof}
We need to show that $\mathrm{WF}(\mathcal{M})$ satisfies all of the axioms of $\mathrm{KP}$.    By Lemma \ref{Th:StandardPartClosedUnderRank}, we are left to verify that $\mathrm{WF}(\mathcal{M})$ satisfies $\Delta_0$-Collection. Let $\phi(x, y, \vec{z})$ be a $\Delta_0$-formula and let $a, \vec{b} \in \mathrm{WF}(M)$ be such that
$$\mathrm{WF}(\mathcal{M}) \models (\forall x \in a) \exists y \phi(x, y, \vec{b}).$$
Since $\mathrm{WF}(\mathcal{M}) \prec_{\Delta_0^\mathcal{P}} \mathcal{M}$,
$$\mathcal{M} \models (\forall x \in a) \exists y \phi(x, y, \vec{b}).$$
Consider $\theta(\gamma, \vec{z}, w)$ defined by
$$(\forall x \in w)(\exists \alpha \in \gamma)\exists y (\phi(x, y, \vec{z})\land \rho(y)= \alpha).$$
Recall that $\Sigma_1$-Collection in $\mathcal{M}$ implies that $\theta(\gamma, \vec{z}, w)$ is equivalent to a $\Sigma_1$-formula. Therefore, using $\Sigma_1$-Foundation, let $\delta \in M$ be the least element of
$$A= \{ \beta \in M \mid \mathcal{M} \models \theta(\beta, \vec{b}, a)\}.$$
Now, every nonstandard $\mathcal{M}$-ordinal is an element of $A$ and so $\delta \in \mathrm{WF}(M)$. Let $\psi(x, z, \gamma, \vec{w})$ be the $\Sigma_1$-formula
$$\exists y \exists \alpha(z= \langle y, \alpha\rangle \land (\alpha \in \gamma)\land \rho(y)= \alpha \land \phi(x, y, \vec{w})).$$
Therefore,
$$\mathcal{M} \models (\forall x \in a) \exists z \psi(x, z, \delta, \vec{b}).$$
Work inside $\mathcal{M}$. By $\Sigma_1$-Collection, there exists $d$ such that $(\forall x \in a) (\exists z\in d) \psi(x, z, \delta, \vec{b})$. Let $c= \mathrm{dom}(d)$. Note that by Lemma \ref{Th:StandardPartClosedUnderRank}, $c^* \subseteq \mathrm{WF}(M)$ and so $c \in \mathrm{WF}(M)$. Hence
$$\mathrm{WF}(\mathcal{M}) \models (\forall x \in a) (\exists y \in c) \phi(x, y, \vec{b}),$$
and therefore $\Delta_0$-Collection holds in $\mathrm{WF}(\mathcal{M})$. So, the Mostowski collapse of $\mathrm{WF}(\mathcal{M})$ witnesses the fact that $\mathrm{WF}(\mathcal{M})$ is isomorphic to an admissible set.
\Square
\end{proof}

In Definitions \ref{Df:CUnbounded} and \ref{contained-def}, we introduce two relationships between models of set theory and their well-founded parts that will be shown to be linked to the existence of proper initial self-embeddings in Sections \ref{Sec:ModelZFCMinusWithNoSelfEmbedding} and \ref{Sec:ModelsWithSelfEmbeddings}.

\begin{Definitions1} \label{Df:CUnbounded}
Let $\mathcal{M} \models \mathrm{KPI}$.

(a) The well-founded part of $\mathcal{M}$ is \textbf{c-bounded} in $\mathcal{M}$, where ``c" stands for ``cardinalitywise", if there is some $x \in M$ such that for all $w \in \mathrm{WF}(M)$  $\mathcal{M} \models |x| > |w|$.

(b) The well-founded part of $\mathcal{M}$ is \textbf{c-unbounded} in $\mathcal{M}$ if the well-founded part of $\mathcal{M}$ is not c-bounded in $\mathcal{M}$, i.e.,
if for all $x \in M$, there exists $w\in \mathrm{WF}(M)$ such that $\mathcal{M} \models |x| \leq |w|$.
\end{Definitions1}

In Section \ref{Sec:ModelZFCMinusWithNoSelfEmbedding} we will see that the well-founded part being c-bounded prevents a model of $\mathrm{KPI}$ from admitting a proper initial self-embedding. In contrast, the following condition will be used in Section \ref{Sec:ModelsWithSelfEmbeddings} to show that nonstandard models of certain extensions of $\mathrm{KPI}$ are guaranteed to admit proper initial self-embedding.

\begin{Definitions1} \label{contained-def}
Let $\mathcal{M} \models \mathrm{KPI}$. We say that the \textbf{well-founded part of $\mathcal{M}$ is contained} if there exists $c \in M$ such that $\mathrm{WF}(M) \subseteq  c^*$.
\end{Definitions1}

The next result shows that if the well-founded part of an $\omega$-standard model is contained, then Theorem \ref{Th:WellFoundedPartIsAdmissible} can be extended to show that the well-founded part satisfies all of the axioms of $\mathrm{KP}^\mathcal{P}$.

\begin{Theorems1} \label{Th:ContainedStandardPartsSatisfyKPP}
Let $\mathcal{M} \models \mathrm{KPI}$ be $\omega$-standard. If the well-founded part of $\mathcal{M}$ is contained, then $\mathrm{WF}(\mathcal{M})\models \mathrm{KP}^\mathcal{P}$.
\end{Theorems1}

\begin{proof}
Suppose that the well-founded part of $\mathcal{M}$ is contained. Let $c \in M$ be such that $\mathrm{WF}(M) \subseteq c^*$. Since $\omega \in \mathrm{WF}(M)$, $\mathrm{WF}(\mathcal{M})$ satisfies the axiom of infinity. By Lemma \ref{Th:StandardPartClosedUnderRank}, we are left to verify that Powerset and $\Delta_0^{\mathcal{P}}$-Collection hold in $\mathrm{WF}(\mathcal{M})$. To see that Powerset holds, let $x \in \mathrm{WF}(M)$. It follows from Lemma \ref{Th:StandardPartClosedUnderRank} that if $y \in M$ with $\mathcal{M} \models y \subseteq x$, then $y \in \mathrm{WF}(M)$. Note that $\Delta_0$-Separation in $\mathcal{M}$ ensures that
$$A= \{ y \in c \mid y \subseteq x\}.$$
is a set. Moreover, $\mathcal{M} \models (A= \mathcal{P}(x))$. Now, $A^* \subseteq \mathrm{WF}(M)$ and so, by Lemma \ref{Th:WellfoundedPartTopless}, $A \in \mathrm{WF}(M)$. Therefore, $\mathrm{WF}(\mathcal{M})\models \textrm{Powerset}$. We are left to verify $\Delta_0^\mathcal{P}$-Collection. Let $\phi(x, y, \vec{z})$ be a $\Delta_0^\mathcal{P}$-formula. Let $a, \vec{b} \in \mathrm{WF}(M)$ be such that
$$\mathrm{WF}(\mathcal{M}) \models (\forall x \in a) \exists y \phi(x, y, \vec{b}).$$
Since $\mathrm{WF}(\mathcal{M}) \subseteq_{\mathrm{topless}}^\mathcal{P} \mathcal{M}$,
$$\mathcal{M}\models (\forall x \in a) (\exists y \in c) \phi(x, y, \vec{b}).$$
Let $\phi^c(x, y, \vec{z})$ be the $\Delta_0$-formula obtained by restricting all of the quantifiers in $\phi(x, y, \vec{z})$ to elements of $c$. Since $\mathrm{WF}(M) \subseteq c^*$ and $\mathrm{WF}(\mathcal{M}) \subseteq_e^\mathcal{P} \mathcal{M}$, for all $x, y, \vec{z} \in \mathrm{WF}(M)$,
$$\mathcal{M} \models \phi(x, y, \vec{z}) \textrm{ if and only if } \mathrm{WF}(\mathcal{M}) \models \phi^c(x, y, \vec{z}).$$
Work inside $\mathcal{M}$. Let
$$B= \{y \in c \mid (\exists x \in a)(\phi^c(x, y, \vec{b}) \land (\forall w \in c)(\phi^c(x, w, \vec{b}) \Rightarrow \rho(y) \leq \rho(w)))\}.$$
Now, $\Delta_1$-Separation ensures that $B$ is a set. Lemma \ref{Th:StandardPartClosedUnderRank} implies that $B^* \subseteq \mathrm{WF}(M)$. Therefore, since $\mathrm{WF}(\mathcal{M}) \subseteq_{\mathrm{topless}}^\mathcal{P} \mathcal{M}$, $B \in \mathrm{WF}(M)$. This shows that
$$\mathrm{WF}(\mathcal{M}) \models (\forall x \in a)(\exists y \in B) \phi(x, y, \vec{b}),$$
and $\mathrm{WF}(\mathcal{M})\models \mathrm{KP}^{\mathcal{P}}$.
\Square
\end{proof}

Note that if $\mathcal{M}$ is a model of $\mathrm{KPI}$ that is not $\omega$-standard, then $\mathrm{WF}(\mathcal{M})$ is isomorphic to $V_\omega$ (the hereditarily finite sets of the metatheory). This observation combined with Theorem \ref{Th:ContainedStandardPartsSatisfyKPP} and Corollary \ref{Th:RanksInKPP} make it clear that if the well-founded part of a model of $\mathrm{KPI}$ is contained, then the well-founded part of that model has access to sequences enumerating the $V_\alpha$s (as in Definition \ref{V_alpha}), thereby yielding the corollary below.



\begin{Coroll1} \label{Th:StandardPartHasRanks} Let $\mathcal{M} \models \mathrm{KPI}$. If the well-founded part of $\mathcal{M}$ is contained, then the sentence expressing
``$V_{\alpha}$ exists for all ordinal $\alpha$'' holds in $\mathrm{WF}(\mathcal{M})$. \end{Coroll1}

We can now see that if the well-founded part of a model of $\mathrm{KPI}$ is contained, then the well-founded part is c-bounded in that model.

\begin{Theorems1} \label{Th:ContainedImpliesBoundedOverKPI}
Let $\mathcal{M} \models \mathrm{KPI}$. If the well-founded part of $\mathcal{M}$ is contained, then the well-founded part of $\mathcal{M}$ is c-bounded in $\mathcal{M}$.
\end{Theorems1}

\begin{proof}
Assume that the well-founded part of $\mathcal{M}$ is contained. Let $C \in M$ be such that $\mathrm{WF}(M) \subseteq C^*$. Suppose, for a contradiction, that the well-founded part of $\mathcal{M}$ is c-unbounded in $\mathcal{M}$. Let $X \in \mathrm{WF}(\mathrm{M})$ be such that $\mathcal{M} \models |C| \leq |X|$. It follows from Theorem \ref{Th:ContainedStandardPartsSatisfyKPP} that there exists $Y \in \mathrm{WF}(M)$ such that $\mathcal{M} \models Y= \mathcal{P}(X)$. Work inside $\mathcal{M}$. By Cantor's Theorem, $|X| < |Y|$. Note the usual proof of Cantor's theorem can be carried out in $\mathrm{KP}$ since it uses $\Delta_1$-Separation, which is a theorem of $\mathrm{KP}$. But, $|Y| \leq |C| \leq |X|$, which is a contradiction. Therefore the well-founded part of $\mathcal{M}$ is c-bounded in $\mathcal{M}$.
\Square
\end{proof}

The next lemma shows that in the special case when a countable model $\mathcal{M}$ of  $\mathrm{KPI}$ is $\omega$-nonstandard, the axioms of $\mathrm{KPI}$ are sufficient to ensure that the well-founded part is contained (note that the well-founded part of an $\omega$-nonstandard  is isomorphic to the hereditarily finite sets).

\begin{Lemma1} \label{Th:OmegaNonstandardModelsHaveContainedStandardPart}
If $\mathcal{M} \models \mathrm{KPI}$ is $\omega$-nonstandard, then the well-founded part of $\mathcal{M}$ is contained.
\end{Lemma1}

\begin{proof}
If $\mathcal{M} \models \mathrm{KPI}$ is $\omega$-nonstandard, then $\mathrm{WF}(\mathcal{M})$ is isomorphic to $V_\omega$ and $\mathrm{WF}(M) \subseteq (L_\omega^\mathcal{M})^*$. \Square
\end{proof}

Similarly, if $\mathcal{M}$ in nonstandard and satisfies all of the axioms of $\mathrm{KP}^\mathcal{P}$, then the well-founded part of $\mathcal{M}$ is contained.

\begin{Lemma1} \label{Th:WellFoundedPartsOfKPPContained}
If $\mathcal{M} \models \mathrm{KP}^\mathcal{P}$ is nonstandard, then the well-founded part of $\mathcal{M}$ is contained.
\end{Lemma1}

\begin{proof}
Let $\mathcal{M}= \langle M, \in^\mathcal{M} \rangle$ be a nonstandard model of $\mathrm{KP}^\mathcal{P}$. Let $x \in M$ be such that $x \notin \mathrm{WF}(M)$. Let $\alpha= \rho^\mathcal{M}(x)$. By Lemma \ref{Th:StandardPartClosedUnderRank}, $\alpha \in M \backslash \mathrm{WF}(M)$ and $\mathrm{WF}(M) \subseteq (V_\alpha^\mathcal{M})^*$.\Square
\end{proof}

Over the theory $\mathrm{ZF}^-+\mathrm{WO}+\forall \alpha\ \Pi_\infty^1-\mathrm{DC}_\alpha$, we get a converse to Theorem \ref{Th:ContainedImpliesBoundedOverKPI}.

\begin{Lemma1} \label{Th:CBoundedImpliesPowersetInWellFoundedPart}
Let $\mathcal{M}$ be a model of $\mathrm{ZF}^-+\mathrm{WO}+\forall \alpha\ \Pi_\infty^1-\mathrm{DC}_\alpha$. If the well-founded part of $\mathcal{M}$ is c-bounded in $\mathcal{M}$, then
$$\mathrm{WF}(\mathcal{M}) \models \mathrm{Powerset}.$$
\end{Lemma1}


\begin{proof}
Suppose that the well-founded part of $\mathcal{M}$ is c-bounded in $\mathcal{M}$. Using the axiom WO and transitive collapse inside $\mathcal{M}$, let $\kappa \in M$ be such that $\mathcal{M} \models (\kappa \textrm{ is a cardinal})$ and for all $X \in \mathrm{WF}(M)$, $\mathcal{M} \models |X| < \kappa$. Since $\mathrm{WF}(\mathcal{M}) \subseteq_{\mathrm{topless}}^{\mathcal{P}} \mathcal{M}$, to see that $\mathrm{WF}(\mathcal{M}) \models \mathrm{Powerset}$ it is sufficient to show that for all $X \in \mathrm{WF}(M)$, $\mathcal{M}$ thinks that $\mathcal{P}(X)$ exists. Suppose, for a contradiction, that $Y \in \mathrm{WF}(M)$ is such that $\mathcal{M}$ believes that the powerset of $Y$ does not exist. Note that if $Z \in M$ with $\mathcal{M} \models Z \subseteq Y$, then $Z \in \mathrm{WF}(M)$. Consider the formula $\phi(f, y, Y, \kappa)$ defined by
$$(\exists \alpha \in \kappa)(f \textrm{ is a function with domain } \alpha) \Rightarrow (y \subseteq Y) \land (y \notin \mathrm{rng}(f)).$$
Now, suppose that there exists $f \in M$ such that $\mathcal{M} \models \forall y \neg \phi(f, y, Y, \kappa)$. It follows that
$$\mathcal{M} \models (f \textrm{ is a function}) \land (\mathrm{rng}(f)= \mathcal{P}(Y)),$$
which contradicts the fact that the powerset of $Y$ does not exist. Therefore
$$\mathcal{M} \models \forall f \exists y \phi(f, y, Y, \kappa),$$
and so, by $\Pi^1_\infty\mathrm{-DC}_\kappa$, there exists $f \in M$ such that
$$\mathcal{M} \models (f \textrm{ is a function}) \land (\mathrm{dom}(f)= \kappa) \land (\forall \alpha \in \kappa)\phi(f\upharpoonright \alpha, f(\alpha), Y, \kappa).$$
Work inside $\mathcal{M}$. If $\alpha \in \kappa$ is least such that $f(\alpha) \nsubseteq Y$, then $\phi(f \upharpoonright \alpha, f(\alpha), Y, \kappa)$ does not hold. Therefore for all $y \in \mathrm{rng}(f)$, $y \subseteq Y$. If $\alpha \in \kappa$ is least such that there exists $\beta \in \alpha$ with $f(\alpha) = f(\beta)$, then $\phi(f \upharpoonright \alpha, f(\alpha), Y, \kappa)$ does not hold. This shows that $f$ is injective. Therefore, we have $\mathrm{rng}(f)^* \subseteq \mathrm{WF}(M)$ and $\mathcal{M} \models \kappa \leq |\mathrm{rng}(f)|$. And, since $\mathrm{WF}(\mathcal{M}) \subseteq_{\mathrm{topless}}^{\mathcal{P}} \mathcal{M}$, $\mathrm{rng}(f) \in \mathrm{WF}(M)$, and the fact that $\mathcal{M} \models \kappa \leq |\mathrm{rng}(f)|$ contradicts our choice of $\kappa$. This shows that $\mathrm{WF}(\mathcal{M})$ satisfies the powerset axiom.
\end{proof}

\begin{Lemma1} \label{Th:PowersetInWellFoundedPartImpliesContained}
Let $\mathcal{M}$ be a model of $\mathrm{ZF}^-+\mathrm{WO}$. If $\mathcal{M}$ is nonstandard and
$$\mathrm{WF}(\mathcal{M}) \models \mathrm{Powerset},$$
then the well-founded part of $\mathcal{M}$ is contained.
\end{Lemma1}

\begin{proof}
Suppose that $\mathcal{M}$ is nonstandard and
$$\mathrm{WF}(\mathcal{M}) \models \mathrm{Powerset}.$$
Consider the formula $\phi(f, \alpha)$ that expresses that $f$ is a function with domain $\alpha$ such that $f(\beta)=V_\beta$ for all $\beta<\alpha$.
Suppose that the class $\{ \alpha \in \mathrm{Ord}^\mathcal{M} \mid \mathcal{M} \models \neg \exists f \phi(f, \alpha)\}$ is nonempty and, using foundation, let $\xi \in \mathrm{Ord}^\mathcal{M}$ be the least element of this class. We claim that $\xi \notin \mathrm{o}(M)$. Suppose, for a contradiction, that $\xi \in \mathrm{o}(M)$.
Since $\mathrm{WF}(\mathcal{M}) \models \mathrm{Powerset}$, $\xi$ is not a successor ordinal and therefore must be a limit ordinal. Work inside $\mathcal{M}$. Using Collection and Separation in $\mathcal{M}$, the class
$$A= \{ f \mid (\exists \alpha \in \xi) \phi(f, \alpha)\}$$
is a set. Now, $\bigcup A$ is a function that satisfies $\phi\left(\bigcup A, \xi\right)$, which contradicts our choice of $\xi$. This shows that if $\{ \alpha \in \mathrm{Ord}^\mathcal{M} \mid \mathcal{M} \models \neg \exists f \phi(f, \alpha)\}$ is nonempty, then its least element cannot be standard. Therefore, since $\mathcal{M}$ is nonstandard, there exists $f \in M$ and $\gamma \in \mathrm{Ord}^\mathcal{M} \backslash \mathrm{o}(M)$ such that $\mathcal{M} \models \phi(f, \gamma)$. Moreover, since $\mathrm{WF}(\mathcal{M}) \subseteq_{\mathrm{topless}}^{\mathcal{P}} \mathcal{M}$, there exists $\nu \in \gamma^*\backslash \mathrm{o}(M)$. A standard induction argument inside $\mathcal{M}$ shows that
$$\mathcal{M} \models \forall x( x \in f(\nu) \iff \rho(x)< \nu).$$
Therefore, by Lemma \ref{Th:StandardPartClosedUnderRank} and the fact that $\nu$ is nonstandard, $\mathrm{WF}(M) \subseteq f(\nu)^*$. This shows that the well-founded part of $\mathcal{M}$ is contained.
\Square
\end{proof}

\begin{Theorems1} \label{Th:CBoundedImpliesContainedOverDependentChoices}
Let $\mathcal{M}$ be a model of $\mathrm{ZF}^-+\mathrm{WO}+\forall \alpha\ \Pi_\infty^1-\mathrm{DC}_\alpha$. If the well-founded part of $\mathcal{M}$ is c-bounded in $\mathcal{M}$, then the well-founded part of $\mathcal{M}$ is contained.
\end{Theorems1}

\begin{proof}
Suppose that the well-founded part of $\mathcal{M}$ is c-bounded in $\mathcal{M}$. Therefore $\mathcal{M}$ is nonstandard and, by Lemma \ref{Th:CBoundedImpliesPowersetInWellFoundedPart},
$$\mathrm{WF}(\mathcal{M})\models \mathrm{Powerset}.$$
The fact that the well-founded part of $\mathcal{M}$ is contained now follows from Lemma \ref{Th:PowersetInWellFoundedPartImpliesContained}.
\Square
\end{proof}

\section[Obstructing initial self-embeddings]{Obstructing initial self-embeddings} \label{Sec:ModelZFCMinusWithNoSelfEmbedding}

In this section we establish the first main result of the paper (Theorems \ref{Th:ModelOfZFCminusWithoutSelfEmbedding}) on the existence of countable nonstandard models of $\mathrm{ZF}^-+\mathrm{WO}$ with no nontrivial initial self-embeddings.  Furthermore, in Theorem \ref{Th:Uncountablemodelswithnoselfembeddings} we
exhibit nonstandard \textit{uncountable} models of ZF with no proper initial self-embeddings.

We begin with verifying that an initial embedding of a model of $\mathrm{KPI}$ must fix the well-founded part of this model. This result will allow us to show that models of $\mathrm{KPI}$ in which the well-founded part is c-unbounded (in the sense of Definition \ref{Df:CUnbounded}) do not admit proper initial self-embedding.

\begin{Lemma1} \label{Th:InitialSelfEmbeddingsIdentityOnStandardPart}
Let $\mathcal{M} \models \mathrm{KPI}$. If $j: \mathcal{M} \longrightarrow \mathcal{M}$ is an initial self-embedding, then $j$ is the identity on $\mathrm{WF}(\mathcal{M})$.
\end{Lemma1}

\begin{proof}
Let $j: \mathcal{M} \longrightarrow \mathcal{M}$ be a proper initial self-embedding. Suppose that $j$ is not the identity on $\mathrm{WF}(\mathcal{M})$ and let $x \in\mathrm{WF}(\mathcal{M})$ be $\in^\mathcal{M}$-least such that $j(x)\neq x$. Now, if $z \in M$ with $\mathcal{M} \models z \in x$ and $\mathcal{M} \models z \notin j(x)$, then $j(z) = z$ and $\mathcal{M} \models (z \in x) \land (j(z) \notin j(x))$, which is a contradiction. Similarly, if $z \in M$ with $\mathcal{M} \models z \notin x$ and $\mathcal{M} \models z \in j(x)$, then $j^{-1}(z)\neq z$ and $\mathcal{M} \models j^{-1}(z) \in x$, which contradicts the fact that $x$ is the $\in^{\mathcal{M}}$-least thing moved by $j$.
\Square
\end{proof}

\begin{Coroll1}
Let $\mathcal{M} \models \mathrm{KPI}$. If $j: \mathcal{M} \longrightarrow \mathcal{M}$ is a proper initial self-embedding, then
$$\mathrm{WF}(\mathcal{M}) \subseteq_{\mathrm{topless}}^\mathcal{P} j[\mathcal{M}].$$
\end{Coroll1}

The next lemma shows that in addition to containing the well-founded part, the fixed point set $\mathrm{Fix}(j)$ of a proper initial self-embedding $j: \mathcal{M} \longrightarrow \mathcal{M}$ also contains all points that are $\Sigma_1$-definable in $\mathcal{M}$ from points in $\mathrm{Fix}(j)$.

\begin{Lemma1} \label{Th:Sigma1DefinablePointsAreFixed}
Let $\mathcal{M} \models \mathrm{KPI}$. Let $j: \mathcal{M} \longrightarrow \mathcal{M}$ be an initial self-embedding. If $x \in M$ is definable in $\mathcal{M}$ by a $\Sigma_1$-formula with parameters from $\mathrm{Fix}(j)$, then $x \in \mathrm{Fix}(j)$.
\end{Lemma1}

\begin{proof}
Suppose that $\phi(z, \vec{y})$ is a $\Sigma_1$-formula, $\vec{a} \in \mathrm{Fix}(j)$ and $x \in M$ is the unique element of $M$ such that
$$\mathcal{M} \models \phi(x, \vec{a}).$$
Therefore, since $\vec{a} \in \mathrm{Fix}(j)$,
$$j[\mathcal{M}] \models \phi(j(x), \vec{a}).$$
And, since $j[\mathcal{M}] \subseteq_e \mathcal{M}$ and $\phi$ is a $\Sigma_1$-formula,
$$\mathcal{M} \models \phi(j(x), \vec{a}).$$
By the uniqueness of $x$, $j(x)= x$ and $x \in \mathrm{Fix}(j)$.
\Square
\end{proof}

In particular, if $j: \mathcal{M} \longrightarrow \mathcal{M}$ is a proper initial self-embedding and $x \in M$ is a point that is $\Sigma_1$-definable in $\mathcal{M}$ from points in the well-founded part of $\mathcal{M}$, then $x$ must be fixed by $j$. This observation allows us to show if the well-founded part of a nonstandard model of $\mathrm{KPI}$ is c-unbounded, then that model admits no proper initial self-embedding.


\begin{Theorems1} \label{Th:NoSelfEmbeddingWhenStandardPartDense}
Let $\mathcal{M} \models \mathrm{KPI}$. If the well-founded part of $\mathcal{M}$ is c-unbounded in $\mathcal{M}$, and $j:\mathcal{M} \longrightarrow \mathcal{M}$ is an initial self-embedding, then $j$ is the identity embedding.
\end{Theorems1}

\begin{proof}
Assume that the well-founded part of $\mathcal{M}$ is c-unbounded in $\mathcal{M}$ and suppose that $j: \mathcal{M} \longrightarrow \mathcal{M}$ is an initial self-embedding. We will show that $j$ must be the identity function. Let $x \in M$. Let $y \in M$ be such that
$$\mathcal{M} \models y= \mathrm{TC}^\mathcal{M}(\{x\}).$$
Let $X \in \mathrm{WF}(M)$ be such that
$$\mathcal{M} \models |y| \leq |X|.$$
Work inside $\mathcal{M}$. Let $f: y \longrightarrow X$ be injective. Let $X^\prime= \mathrm{rng}(f)$. Consider
$$Y= \{\langle u, v\rangle \in X \times X \mid (u, v \in \mathrm{rng}(f))\land (f^{-1}(u) \in f^{-1}(v))\}.$$
Now, $Y, X^\prime \in \mathrm{WF}(M)$. Consider $\phi(z, W, Z)$ defined by
$$\exists w\exists f \left(\begin{array}{c}
(f:w \longrightarrow W \textrm{ is a bijection})\land  (\forall u, v \in w)(u \in v \iff \langle f(u), f(v)\rangle \in Z)\\
\land (Z\subseteq W \times W)\land \bigcup w \subseteq w \land (z \in w)\land (\forall u \in w)(z \notin u)
\end{array}\right).$$
Note that $\phi(z, W, Z)$ is a $\Sigma_1$-formula and $x$ is the unique point in $M$ such that
$$\mathcal{M} \models \phi(x, X^\prime, Y).$$
Therefore, by Lemmas \ref{Th:Sigma1DefinablePointsAreFixed} and \ref{Th:Sigma1DefinablePointsAreFixed}, $j(x)=x$. Since $x \in M$ was arbitrary, this shows that $j$ is the identity embedding, as desired.
\Square
\end{proof}

This allows us to show that there are nonstandard $\omega$-standard models of $\mathrm{ZF}^-+\mathrm{WO}$ that are not isomorphic to a transitive subclass of themselves. To build such a model we will employ the following consequence of \cite[Theorem 2.2]{fri73}.

\begin{Lemma1} \label{Th:NonstandardModelWithSameStandardOrdinals}
Let $T$ be a recursive $\mathcal{L}$-theory. If $A$ is a countable admissible set such that $\langle A, \in \rangle \models T$, then there exists a nonstandard $\mathcal{L}$-structure $\mathcal{M}$ such that $\mathcal{M} \models \mathrm{KP}+ T$ and $\mathrm{o}(M)= A\cap \mathrm{Ord}$.
\end{Lemma1}

\begin{Theorems1} \label{Th:ModelOfZFCminusWithoutSelfEmbedding}
There exists a countable nonstandard $\omega$-standard model $\mathcal{M} \models \mathrm{ZF}^-+\mathrm{WO}$ such that there is no initial self-embedding $j:\mathcal{M} \longrightarrow \mathcal{M}$ other than the identity embedding.
\end{Theorems1}

\begin{proof}
Note that $\langle H_{\aleph_1}, \in \rangle \models \mathrm{ZF}^-+\mathrm{WO}+\forall x (|x| \leq \aleph_0)$. Therefore, by the Downwards L\"{o}wenheim-Skolem Theorem and the Mostowski Collapsing Lemma, there exists a countable admissible set $A$ such that $\omega \in A$ and $\langle A, \in \rangle \equiv \langle H_{\aleph_1}, \in \rangle$. So, by Lemma \ref{Th:NonstandardModelWithSameStandardOrdinals} and the Downwards L\"{o}wenheim-Skolem Theorem, there exists a countable $\mathcal{L}$-structure $\mathcal{M}$ such that $\mathcal{M} \models \mathrm{ZF}^-+\mathrm{WO}+\forall x (|x| \leq \aleph_0)$, $\mathcal{M}$ is nonstandard and $\mathrm{o}(M)= A \cap \mathrm{Ord}$. Since $\omega \in A$, the well-founded part of $\mathcal{M}$ is c-unbounded in $\mathcal{M}$ and so, by Theorem \ref{Th:NoSelfEmbeddingWhenStandardPartDense}, there is no proper initial self-embedding $j:\mathcal{M} \longrightarrow \mathcal{M}$.
\Square
\end{proof}

We conclude this section by exhibiting uncountable nonstandard models of $\mathrm{ZF}$ that carry no proper initial self-embeddings. Before doing so, let us note that it is well-known that every consistent extension of $\mathrm{ZF}$ has a model of cardinality $\aleph_1$ that carries no proper \textit{rank-initial} self-embedding. To see this, recall that by a classical result due to Keisler and Morely (first established in \cite{Keisler-Morley}, and exposited as Theorem 2.2.18 of \cite{Chang-Keisler}) every countable model of $\mathrm{ZF}$ has a proper elementary end-extension.  It is easy to see that an elementary end-extension of a model of $\mathrm{ZF}$ is a rank-extension. Now if $T$ is a consistent extension of $\mathrm{ZF}$, we can readily build a \textit{countable nonstandard} model of $T$ and use the Keisler-Morley theorem $\aleph_1$-times (while taking unions at limit ordinals) to build a so-called $\aleph_1$-like model of $T$, i.e., a model $\cal{M}$ of power $\aleph_1$ such that $a^*$ is finite or countable for each $a \in M$. It is evident that $\cal{M}$ is nonstandard. Moreover, $\cal{M}$ carries no proper rank initial embedding $j$ since any such embedding $j$ would have to have the property that $j[\mathcal{M}]$ is a submodel of some structure of the form $V^{\mathcal{M}}_{\alpha}$ for some ``ordinal" $\alpha$ of $\mathcal{M}$, which is impossible, since $(V^{\mathcal{M}}_{\alpha})^*$ is countable thanks to the fact that $\cal{M}$ is $\aleph_1$-like.

\begin{Theorems1} \label{Th:Uncountablemodelswithnoselfembeddings} Every consistent extension of $\mathrm{ZF}+ V = L$ has a nonstandard model of power $\aleph_1$ that carries no proper initial self-embedding.
\end{Theorems1}

\begin{proof} Let $T$ be a consistent extension of $\mathrm{ZF} + V = L$, and $\mathcal{M}$ be a nonstandard $\aleph_1$-like model of $T$. Recall that, provably in $\mathrm{ZF} + V = L$, there is a $\Sigma_1$-formula $\sigma(x,y)$ that describes the graph of a bijection $f$ between the class $V$ of sets and the class $\mathrm{Ord}$ of ordinals; see, e.g., the proof of Lemma 13.19 of \cite{Jechbook}. Suppose $j$ is an initial embedding of $\mathcal{M}$. We will show that $j$ is not a proper embedding by verifying that every element $m$ of $\mathcal{M}$ is in the image of $j$. Being an ordinal is a $\Delta_0$-property and thus preserved by $j$, so $j[\mathrm{Ord}^\mathcal{M}] \subseteq \mathrm{Ord}^\mathcal{M}$. Since $\mathcal{M}$ is $\aleph_1$-like, $j[\mathrm{Ord}^\mathcal{M}]$ must be cofinal in $\mathrm{Ord}^\mathcal{M}$. But, $j$ is initial, so $j[\mathrm{Ord}^\mathcal{M}]= \mathrm{Ord}^\mathcal{M}$.
To show that $j[\mathcal{M}] = \mathcal{M}$, suppose $m \in M$. Then there is a unique $\alpha \in \mathrm{Ord}^{\mathcal{M}}$ such that $\mathcal{M} \models \sigma(m, \alpha)$. Since every ordinal of $\mathcal{M}$ is in $j[\mathrm{Ord}^{\mathcal{M}}]$, there is some $\beta \in \mathrm{Ord}^{\mathcal{M}}$ such that $j(\beta) = \alpha$. Let $m_0$ be the unique element of $\mathcal{M}$ such that  $\mathcal{M} \models \sigma(m_0, \beta)$. Then since $j$ is an embedding,  $j[\mathcal{M}] \models \sigma(j(m_0), j(\beta))$, and by the choice of $\beta$, $j[\mathcal{M}] \models \sigma(j(m_0), \alpha))$, which coupled with the fact that $\sigma$ is a $\Sigma_1$-formula, yields $\mathcal{M} \models \sigma(j(m_0), \alpha)$. So in light of the fact that $\sigma$ within $\mathcal{M}$ defines the graph of a bijection $f$ between $V$ and $\mathrm{Ord}$, and $f(m)=\alpha$ (by the choice of $\alpha$), we can can conclude that $j(m_0)=m$, thereby showing that $j[\mathcal{M}] = \mathcal{M}$. \Square
\end{proof}

\section[Constructing initial self-embeddings]{Constructing initial self-embeddings} \label{Sec:ModelsWithSelfEmbeddings}

In the previous section we saw that if the well-founded part of a model $\mathcal{M}$ of $\mathrm{KPI}$ is c-unbounded (in the sense of Definition \ref{Df:CUnbounded}) in $\mathcal{M}$, then there is no proper initial self-embedding of $\mathcal{M}$. In this section we prove an adaption of H. Friedman's Self-embedding Theorem \cite[Theorem 4.1]{fri73} that ensures the existence of proper initial self-embeddings of models of extensions of $\mathrm{KPI}$ with contained well-founded parts.\footnote{It is known that $\mathrm{KP+\lnot I}$ + $\Sigma_1$-separation is bi-interpretable with the fragment $\mathrm{I}\Sigma_1$ of PA (Peano arithmetic), where $\mathrm{KP+\lnot I}$ is KP plus the negation the axiom of infinity; the proof is implicit in the proof of the main result of Kaye and Wong \cite{kw07}. Moreover, the bi-interpretation at work makes it clear that the study of initial self-embeddings of models of $\mathrm{KP+\lnot I}$ + $\Sigma_1$-Separation boils down to the study of initial self-embeddings of models of $\mathrm{I}\Sigma_1$. The interested reader can consult \cite{BahEna18} for a systematic study of initial self-embeddings of $\mathrm{I}\Sigma_1$.}




We now turn to the investigation of conditions under which models of $\mathrm{KPI}$ with contained well-founded parts admit proper initial self-embeddings. We begin with the verification that $\Delta_1$-Separation ensures that $\Sigma_0$-types with parameters from the well-founded part that are realised are coded in the standard system; and that for $n>0$, $\Sigma_n$-Separation is sufficient to ensure the corresponding condition for $\Sigma_n$-types.


\begin{Lemma1} \label{Th:Sigma1TypesWithParametersAreCoded} \footnote{We are grateful to Kameryn Williams for spotting an unnecessary assumption in an earlier version of this Lemma \ref{Th:Sigma1TypesWithParametersAreCoded}}
Suppose $n \in \omega$ and $\mathcal{M}$ is a model of $\mathrm{KPI}+\Sigma_{n}\textrm{-Separation}$ such that the well-founded part of $\mathcal{M}$ is contained. If $\vec{a} \in M$, then
$$\{\langle \ulcorner \phi(x, \vec{y}) \urcorner, b\rangle \mid \phi \textrm{ is } \Sigma_{n}, \ b \in \mathrm{WF}(M) \textrm{ and } \mathcal{M} \models \phi(b, \vec{a}) \} \in \mathrm{SSy}(\mathcal{M}).$$
\end{Lemma1}

\begin{proof}
Let $C \in M$ be such that $\mathrm{WF}(M) \subseteq C^*$. Let $a_1, \ldots, a_l \in M$. Work inside $\mathcal{M}$. Consider
$$D= \{\langle q, b\rangle \in C \mid (q \in \omega) \land (a= \langle b, a_1, \ldots, a_l \rangle) \land \mathrm{Sat}_{\Sigma_n}(q, a) \}.$$
Thanks to Theorem \ref{Complexityofpartialsat}, $\Sigma_n$-Separation ($\Delta_1$-Separation when $n=0$) ensures that $D$ is a set in $\mathcal{M}$. It is clear that $D$ codes
$$\{\langle \ulcorner \phi(x, \vec{y}) \urcorner, b\rangle \mid \phi \textrm{ is } \Sigma_n, b \in \mathrm{WF}(M) \textrm{ and } \mathcal{M} \models \phi(b, \vec{a}) \}.$$
\Square
\end{proof}


As verified in the next lemma, in the special case when the model is not $\omega$-standard, in Lemma \ref{Th:Sigma1TypesWithParametersAreCoded} the assumption that the well-founded part is contained can be dropped and the assumption that $\Sigma_n$-Separation holds can be replaced by a fragment of the collection scheme coupled with a fragment foundation scheme.

\begin{Lemma1} \label{Th:SigmaTypesOfOmegaNonstandardModelsAreCoded}
Suppose that $n \in \omega$, $m= \max\{1, n\}$, and $\mathcal{M}$ is an $\omega$-nonstandard model of $ \mathrm{KPI}+ \Pi_{m-1}\textrm{-Collection}+\Pi_{n+1}\textrm{-Foundation}$. If $\vec{a} \in M$, then
$$\{\langle \ulcorner \phi(x, \vec{y}) \urcorner, b\rangle \mid \phi \textrm{ is } \Sigma_n, b \in \mathrm{WF}(M) \textrm{ and } \mathcal{M} \models \phi(b, \vec{a}) \} \in \mathrm{SSy}(\mathcal{M}).$$
\end{Lemma1}

\begin{proof}
By Lemma \ref{Th:OmegaNonstandardModelsHaveContainedStandardPart}, the well-founded part of $\mathcal{M}$ is contained. Let $C \in M$ be such that $\mathrm{WF}(M) \subseteq C^*$. Let $a_1, \ldots, a_m \in M$.  Let $\psi_1(\alpha, f,C)$ be the formula:

$$\left(\begin{array}{c}
\alpha \in \omega \land \mathrm{dom}(f)= \alpha+1 \land f(0)= \emptyset \ \land \\
(\forall \beta \in \mathrm{dom}(f))(\forall y \in C)(y \in f(\beta+1) \iff y \subseteq f(\beta))

\end{array}\right).$$

\noindent Consider the formula $\theta(\alpha, \omega, a_1, \ldots, a_m)$ defined by:

$$\exists v (\exists f \in C) \left( \psi_1(\alpha, f,C) \land \psi_2(\alpha, f,C,v,a_1,\ldots,a_m)\right),$$ where $\psi_2(\alpha, f,C,v, a_1,\ldots,a_m) $ is:

$$\left(\begin{array}{c}

(\forall \langle q, b\rangle \in C)\left(\begin{array}{c}
\langle q, b \rangle \in v \cap f(\alpha) \iff\\
(\langle q, b\rangle \in f(\alpha))\land \mathrm{Sat}_{\Sigma_n}(q, \langle b, a_1, \ldots, a_m \rangle)
\end{array}\right)
\end{array}\right).$$



\noindent By $\Pi_{m-1}$-Collection, the formula $\theta(\alpha, \omega, a_1, \ldots, a_m)$ is equivalent to a $\Sigma_{n+1}$-formula. Now, since $\mathrm{WF}(\mathcal{M})$ is isomorphic to $V_\omega$, for all $\alpha \in \mathrm{o}(M)$,
$$\mathcal{M} \models \theta(\alpha, \omega, a_1, \ldots, a_m).$$
Therefore, by $\Pi_{n+1}$-Foundation, there exists $\gamma \in (\omega^\mathcal{M})^* \backslash \mathrm{WF}(M)$ such that
$$\mathcal{M} \models \theta(\gamma, \omega, a_1, \ldots, a_m).$$

\noindent Let $v \in M$ be such that
$$\mathcal{M} \models (\exists f \in C) \left( \psi_1(\gamma, f,C) \land \psi_2(\gamma, f,C,v,a_1,\ldots,a_m)\right).$$



\noindent Since $\gamma \in (\omega^\mathcal{M})^*$ is nonstandard, it follows that $v$ codes
$$\{\langle \ulcorner \phi(x, \vec{y}) \urcorner, b\rangle \mid \phi \textrm{ is } \Sigma_n, b \in \mathrm{WF}(M) \textrm{ and } \mathcal{M} \models \phi(b, \vec{a}) \} \in \mathrm{SSy}(\mathcal{M}).$$
\Square
\end{proof}

Lemma \ref{Th:Sigma1TypesWithParametersAreCoded} allows us to prove the following theorem that gives a sufficient condition for nonstandard models of extensions of $\mathrm{KPI}$ to admit proper initial self-embeddings.


\begin{Theorems1} \label{Th:MainSelfEmbeddingResult2}
Let $p \in \omega$, $\mathcal{M}$ be a countable model of $\mathrm{KPI}+\Sigma_{p+1}\textrm{-Separation}+\Pi_p\textrm{-Collection}$, and let $b, B \in M$ and $c \in B^*$ with the following properties:
\begin{itemize}
\item[(I)] $\mathcal{M} \models \bigcup B \subseteq B$.
\item[(II)] $\mathrm{WF}(M) \subseteq B^*$.
\item[(III)] for all $\Pi_p$-formulae $\phi(\vec{x}, y, z)$ and for all $a \in \mathrm{WF}(M)$,
$$\textrm{if } \mathcal{M} \models \exists \vec{x} \phi(\vec{x}, a, b), \textrm{ then } \mathcal{M} \models (\exists \vec{x} \in B)\phi(\vec{x}, a, c).$$
\end{itemize}
Then there exists a proper initial self-embedding $j: \mathcal{M} \longrightarrow \mathcal{M}$ such that $j[\mathcal{M}] \subseteq_e B^*$, $j(b)=c$ and $j[\mathcal{M}] \prec_p \mathcal{M}$.
\end{Theorems1}

\begin{proof}
It follows from (II) that $B \in M$ witnesses the fact that the well-founded part of $\mathcal{M}$ is contained. Let $\langle d_i \mid i \in \omega\rangle$ be an enumeration of $M$ such that $d_0= b$. Let $\langle e_i \mid i \in \omega\rangle$ be an enumeration of $B^*$ in which every element of $B^*$ appears infinitely often. We will construct an initial embedding $j: \mathcal{M} \longrightarrow \mathcal{M}$ by constructing sequences $\langle u_i \mid i \in \omega\rangle$ of elements of $M$ and $\langle v_i \mid i \in \omega \rangle$ of elements of $B^*$ and defining $j(u_i)= v_i$ for all $i \in \omega$. After stage $n \in \omega$, we will have chosen $u_0, \ldots, u_{n} \in M$ and $v_0, \ldots, v_{n} \in B^*$ and maintained
\begin{itemize}
\item[]($\dagger_n$) for all $\Pi_{p}$-formulae, $\phi(\vec{x}, z, y_0, \ldots, y_{n})$, and for all $a \in \mathrm{WF}(M)$,
$$\textrm{if } \mathcal{M} \models \exists \vec{x} \phi(\vec{x}, a, u_0, \ldots, u_{n}), \textrm{ then } \mathcal{M} \models (\exists \vec{x} \in B)\phi(\vec{x}, a, v_0, \ldots, v_{n}).$$
\end{itemize}
At stage $0$, let $u_0= b$ and let $v_0= c$. By (III), this choice of $u_0$ and $v_0$ satisfy ($\dagger_0$). Let $n \in \omega$ with $n \geq 1$. Assume that we have chosen $u_0, \ldots, u_{n-1} \in M$ and $v_0, \ldots, v_{n-1} \in B^*$ and that ($\dagger_{n-1}$) holds.\\
{\bf Case $n=2k+1$ for $k \in \omega$:} This step will ensure that the embedding $j: \mathcal{M} \longrightarrow \mathcal{M}$ is initial. If
$$\mathcal{M}\models e_k \notin \mathrm{TC}(\{v_0, \ldots, v_{n-1}\}),$$
then let $u_n= u_0$ and $v_n= v_0$. This choice of $u_n$ and $v_n$ ensure that $u_0, \ldots, u_n \in M$ and $v_0, \ldots, v_n \in B^*$ satisfy ($\dagger_{n}$). If
$$\mathcal{M}\models e_k \in \mathrm{TC}(\{v_0, \ldots, v_{n-1}\}),$$
then let $v_n= e_k$ and we need to choose $u_n$ to satisfy ($\dagger_{n}$). By Lemma \ref{Th:Sigma1TypesWithParametersAreCoded},
$\Sigma_{p+1}$-separation, and $\Pi_p$-collection, there exists $D \in M$ that codes the class
$$\{ \langle \ulcorner\phi(\vec{x}, z, y_0, \ldots, y_{n}) \urcorner, a \rangle \mid a \in \mathrm{WF}, \phi \textrm{ is } \Sigma_{p} \textrm{ and } \mathcal{M} \models (\forall \vec{x} \in B)\phi(\vec{x}, a, v_0, \ldots, v_n)\} \in \mathrm{SSy}(\mathcal{M}).$$
By Corollary \ref{Th:StandardPartHasRanks}, the well-founded part of $\mathcal{M}$ believes that ranks exist. For all $\alpha \in \mathrm{o}(M)$, let $D_\alpha \in M$ be such that
$$\mathcal{M} \models D_\alpha= D \cap V_\alpha.$$
Note that for all $\alpha \in \mathrm{o}(M)$, $D_\alpha \in \mathrm{WF}(M)$. We have that for all $\alpha \in \mathrm{o}(M)$,
\begin{equation}  \label{eq:BackStepOfEmeddingTheoremClaim3}
\mathcal{M} \models \exists v((v \in \mathrm{TC}(\{v_0, \ldots, v_{n-1}\}) \land (\forall \langle m, a \rangle \in D_\alpha) (\forall \vec{x} \in B) \mathrm{Sat}_{\Sigma_{p}}(m, \langle \vec{x}, a, v_0, \ldots, v_{n-1}, v \rangle)).
\end{equation}
{\bf Claim:} For all $\alpha \in \mathrm{o}(M)$,
\begin{equation} \label{eq:BackStepOfEmeddingTheoremClaim1}
\mathcal{M} \models (\exists v \in \mathrm{TC}(\{u_0, \ldots, u_{n-1}\}))(\forall \langle m, a\rangle \in D_\alpha)\forall \vec{x}\ \mathrm{Sat}_{\Sigma_p}(m, \langle \vec{x}, a, u_0, \ldots, u_{n-1}, v\rangle).
\end{equation}
To prove this claim, suppose not, and let $\alpha \in \mathrm{o}(M)$ be such that
$$\mathcal{M} \models (\forall v \in \mathrm{TC}(\{u_0, \ldots, u_{n-1}\}))(\exists \langle m, a\rangle \in D_\alpha) \exists \vec{x}\ \neg \mathrm{Sat}_{\Sigma_p}(m, \langle \vec{x}, a, u_0, \ldots, u_{n-1}, v \rangle).$$
By $\Pi_p$-collection,
\begin{equation} \label{eq:BackStepOfEmeddingTheoremClaim2}
\mathcal{M} \models \exists C (\forall v \in \mathrm{TC}(\{u_0, \ldots, u_{n-1}\}))(\exists \langle m, a\rangle \in D_\alpha) (\exists \vec{x} \in C)\ \neg \mathrm{Sat}_{\Sigma_p}(m, \langle \vec{x}, a, u_0, \ldots, u_{n-1}, v \rangle).
\end{equation}
Now, $\Pi_p$-Collection implies that (\ref{eq:BackStepOfEmeddingTheoremClaim2}) is equivalent to a $\Sigma_{p+1}$-formula. Therefore, by ($\dagger_{n-1}$),
$$\mathcal{M} \models (\exists C \in B)(\forall v \in \mathrm{TC}(\{v_0, \ldots, v_{n-1}\}))(\exists \langle m, a\rangle \in D_\alpha) (\exists \vec{x} \in C)\ \neg \mathrm{Sat}_{\Sigma_p}(m, \langle \vec{x}, a, v_0, \ldots, v_{n-1}, v \rangle).$$
But then
$$\mathcal{M} \models (\forall v \in \mathrm{TC}(\{v_0, \ldots, v_{n-1}\}))(\exists \langle m, a\rangle \in D_\alpha) (\exists \vec{x} \in B)\ \neg \mathrm{Sat}_{\Sigma_p}(m, \langle \vec{x}, a, v_0, \ldots, v_{n-1}, v \rangle),$$
which contradicts (\ref{eq:BackStepOfEmeddingTheoremClaim3}). This proves the claim.\\



\noindent Consider the formula $\theta(\alpha, D, B, u_0, \ldots, u_{n-1})$ defined by:
$$ (\exists f \in B) \left( \psi_1(\alpha,f,B) \land \psi_2(\alpha,f, D, u_0, \ldots, u_{n-1}) \right),$$
where $\psi_1(\alpha,f,B)$ is:

$$\left(\begin{array}{c}
(\alpha \textrm{ is an ordinal}) \land \mathrm{dom}(f)= \alpha \land f(0)= \emptyset \ \land \\
(\forall \beta \in \mathrm{dom}(f))\left((\beta \textrm{ is a limit ordinal}) \Rightarrow f(\beta)= \bigcup_{\xi < \beta} f(\xi) \right) \land\\
(\forall \beta \in \mathrm{dom}(f))(\forall y \in B)(y \in f(\beta+1) \iff y \subseteq f(\beta))

\end{array}\right),$$
and $\psi_2(\alpha,f, D, u_0, \ldots, u_{n-1})$ is:
$$ (\exists v \in \mathrm{TC}(\{u_0, \ldots, u_{n-1}\}))(\forall \langle m, a\rangle \in D \cap f(\alpha))\forall \vec{x}\ \mathrm{Sat}_{\Sigma_p}(m, \langle \vec{x}, a, u_0, \ldots, u_{n-1}, v\rangle).$$

\noindent Now, for all $\alpha \in \mathrm{o}(M)$,
$$\mathcal{M} \models \theta(\alpha, D, B, u_0, \ldots, u_{n-1}).$$
And $\Pi_p$-Collection implies that $\theta(\alpha, D, B, u_0, \ldots, u_{n-1})$ is equivalent to a $\Pi_{p+1}$-formula. Therefore, by $\Sigma_{p+1}$-Foundation, there exists $\gamma \in \mathrm{Ord}^\mathcal{M} \backslash \mathrm{o}(M)$ and $f \in M$ such that
$$\mathcal{M} \models \psi_1(\gamma,f,B).$$


\noindent and
$$\mathcal{M} \models \psi_2(\gamma, f,D,u_0,\ldots,u_{n-1}).$$

\noindent Let $v \in M$ be such that
$$\mathcal{M} \models v \in \mathrm{TC}(\{u_0, \ldots, u_{n-1}\}),$$
and
$$\mathcal{M} \models (\forall \langle m, a\rangle \in D \cap f(\gamma))\forall \vec{x}\ \mathrm{Sat}_{\Sigma_p}(m, \langle \vec{x}, a, u_0, \ldots, u_{n-1}, v\rangle).$$
Let $u_n= v$. This choice of $u_n$ ensures that $u_0, \ldots, u_n \in M$ and $v_0, \ldots, v_n \in B^*$ satisfy ($\dagger_{n}$).\\
\\
{\bf Case $n=2k$ for $k \in \omega$:} Let $u_n= d_k$. This choice will ensure that the domain of $j$ is all of $M$. By Lemma \ref{Th:Sigma1TypesWithParametersAreCoded}, there exists $A \in M$ that codes the class
$$\{ \langle\ulcorner \phi(\vec{x}, z, y_0, \ldots, y_{n})\urcorner, a\rangle \mid a \in \mathrm{WF}, \phi \textrm{ is } \Pi_{p} \textrm{ and } \mathcal{M} \models \exists \vec{x} \phi(\vec{x}, a, u_0, \ldots, u_{n})\} \in \mathrm{SSy}(\mathcal{M}).$$
Now, by Corollary \ref{Th:StandardPartClosedUnderRank}, the well-founded part of $\mathcal{M}$ believes that ranks exist. For all $\alpha \in  \mathrm{o}(M)$, let $A_\alpha \in M$ be such that
$$\mathcal{M} \models A_\alpha= V_\alpha \cap A.$$
Note that for all $\alpha \in  \mathrm{o}(M)$, $A_\alpha \in \mathrm{WF}(M)$. We have that for all $\alpha \in \mathrm{o}(M)$,
$$\mathcal{M} \models (\forall \langle m, a \rangle \in A_\alpha) \exists \vec{x} \ \mathrm{Sat}_{\Pi_p}(m, \langle \vec{x}, a, u_0, \ldots, u_n\rangle).$$
So, for all $\alpha \in \mathrm{o}(M)$,
$$\mathcal{M}\models \exists v (\forall \langle m, a\rangle \in A_\alpha) \exists \vec{x} \  \mathrm{Sat}_{\Pi_p}(m, \langle \vec{x}, a, u_0, \ldots, u_{n-1}, v \rangle),$$
and, using $\Pi_p$-Collection, this formula is equivalent to a $\Sigma_{p+1}$-formula with parameters $A_\alpha \in \mathrm{WF}(M)$ and $u_0, \ldots, u_{n-1}$. Therefore, by ($\dagger_{n-1}$) and (I), for all $\alpha \in \mathrm{o}(M)$,
$$\mathcal{M} \models (\exists v \in B)(\forall \langle m, a\rangle \in A_\alpha) (\exists \vec{x} \in B) \mathrm{Sat}_{\Pi_p}(m, \langle \vec{x}, a, v_0, \ldots, v_{n-1}, v \rangle).$$

\noindent Consider the formula $\theta(\alpha, A, B, v_0, \ldots, v_{n-1})$ defined by $$(\exists v, f \in B)\left(\psi_1(\alpha,f,B) \land \psi_2(\alpha,A,B,v_0,\ldots,v_{n-1},v) \right),$$ where $\psi_1(\alpha,f,B)$ is as in the proof of Case $(2n=k+1)$, and $\psi_2(\alpha,A,B,v_0,\ldots,v_{n-1},v)$ is:

$$(\forall \langle m, a \rangle \in A \cap f(\alpha))((\exists \vec{x} \in B) \mathrm{Sat}_{\Pi_p}(m, \langle\vec{x}, a, v_0, \ldots, v_{n-1}, v \rangle)) .$$

\noindent Note that $\theta(\alpha, A, B, v_0, \ldots, v_{n-1})$ is equivalent to a $\Pi_{p}$-formula and for all $\alpha \in \mathrm{o}(M)$,
$$\mathcal{M} \models \theta(\alpha, A, B, v_0, \ldots, v_{n-1}).$$
Therefore, by $\Sigma_{p}$-Foundation, there exists $\gamma \in \mathrm{Ord}^\mathcal{M} \backslash \mathrm{o}(M)$ such that:
$$\mathcal{M} \models \theta(\gamma, A, B, v_0, \ldots, v_{n-1}).$$
Let $f, v \in B^*$ be such that

$$\mathcal{M} \models \left( \psi_1(\gamma,f,B^*) \land \psi_2(\gamma, f,A,B,v_0,\ldots,v_{n-1},v) \right),$$

\noindent and let $v_n= v$. Therefore
$$\mathcal{M} \models (\forall \langle m, a \rangle \in A \cap f(\gamma))((\exists \vec{x} \in B) \mathrm{Sat}_{\Pi_p}(m, \langle \vec{x}, a, v_0, \ldots, v_n \rangle)),$$
and this choice of $v_n$ ensures that $u_0, \ldots, u_n \in M$ and $v_0, \ldots, v_n \in B^*$ satisfy ($\dagger_{n}$). This completes the case where $n=2k$ and shows that we can construct sequences $\langle u_i \mid i \in \omega\rangle$ and $\langle v_i \mid i \in \omega \rangle$ while maintaining the conditions ($\dagger_n$) at each stage of the construction. Now, define $j: \mathcal{M} \longrightarrow \mathcal{M}$ by: for all $i \in \omega$, $j(u_i)= v_i$. Our ``back-and-forth" construction ensures that $j$ is a proper initial self-embedding with $j[\mathcal{M}] \subseteq_e B^*$, $j(b)= c$ and $j[\mathcal{M}] \prec_p \mathcal{M}$.
\Square
\end{proof}

In the proof of Theorem \ref{Th:MainSelfEmbeddingResult2}, the only use of $\Sigma_{p+1}$-Separation is to prove $\Sigma_{p+1}$- and $\Pi_{p+1}$-Foundation, and to satisfy the assumptions of Lemma \ref{Th:Sigma1TypesWithParametersAreCoded}. Therefore, in the special case where the model involved is $\omega$-nonstandard, we can replace Lemma \ref{Th:Sigma1TypesWithParametersAreCoded} with Lemma \ref{Th:SigmaTypesOfOmegaNonstandardModelsAreCoded} to obtain the following simplified variant of Theorem \ref{Th:MainSelfEmbeddingResult2}.

\begin{Theorems1} \label{Th:MainSelfEmbeddingResultOmegaNonStandard2}
Let $p \in \omega$, $\mathcal{M}$ be a countable $\omega$-nonstandard model of $\mathrm{KPI}+\Pi_p\textrm{-Collection}+\Pi_{p+2}\textrm{-Foundation}$, and let $b, B \in M$ and $c \in B^*$ with the following properties:
\begin{itemize}
\item[(I)] $\mathcal{M} \models \bigcup B \subseteq B$,
\item[(II)] $\mathrm{WF}(M) \subseteq B^*$, and
\item[(III)] for all $\Pi_p$-formulae $\phi(\vec{x}, z)$,
$$\textrm{if } \mathcal{M} \models \exists \vec{x} \phi(\vec{x}, b), \textrm{ then } \mathcal{M} \models (\exists \vec{x} \in B) \phi(\vec{x}, c).$$
\end{itemize}
Then there exists a proper initial self-embedding $j: \mathcal{M} \longrightarrow \mathcal{M}$ such that $j[\mathcal{M}] \subseteq_e B^*$, $j(b)= c$ and $j[\mathcal{M}] \prec_p \mathcal{M}$.
\end{Theorems1}

Equipped with Theorems \ref{Th:MainSelfEmbeddingResult2} and \ref{Th:MainSelfEmbeddingResultOmegaNonStandard2}, we are now able to demonstrate in Theorem \ref {Th:ExistenceOfElementarySelfEmbeddingsForKPI} and Corollary \ref{Th:ElementarySelfEmbeddingsForZFCminus} that a variety of nonstandard models of $\mathrm{KPI}$ are isomorphic to $\Sigma_n$-elementary transitive substructures of themselves.

\begin{Theorems1} \label{Th:ExistenceOfElementarySelfEmbeddingsForKPI}
Let $p \in \omega$, $\mathcal{M}$ be a countable nonstandard model of $\mathrm{KPI}+\Sigma_{p+1}\textrm{-Separation}+\Pi_p\textrm{-Collection}$ such that the well-founded part of $\mathcal{M}$ is contained, and let $b \in M$. Then there exists a proper initial self-embedding $j: \mathcal{M} \longrightarrow \mathcal{M}$ such that $b \in \mathrm{rng}(j)$ and $j[\mathcal{M}] \prec_p \mathcal{M}$.
\end{Theorems1}

\begin{proof}
Let $C \in M$ be such that $\mathrm{WF}(M) \subseteq C^*$. Work inside $\mathcal{M}$. Consider the formula $\theta(x, y, b)$ defined by
$$x= \langle m, a \rangle \land (m \in \omega) \land y= \langle y_1, \ldots, y_k\rangle \land \mathrm{Sat}_{\Pi_p}(m, \langle y_1, \ldots, y_k, a, b\rangle).$$
Note that if $p \geq 1$, then $\theta(x, y, b)$ is equivalent to a $\Pi_p$-formula, and if $p=0$, then $\theta(x, y, b)$ is equivalent to a $\Sigma_1$-formula. By Lemma \ref{basicimplications}, Strong $\Pi_p$-Collection holds in $\mathcal{M}$. Therefore, there exists a set $D$ such that
$$(\forall x \in C)(\exists y \theta(x, y) \Rightarrow (\exists y \in D)\theta(x, y)).$$
Let $B= \mathrm{TC}(D)$. Now, $B \in M$ is transitive set in $\mathcal{M}$ with $\mathrm{WF}(M) \subseteq B^*$, and for all $\Pi_p$-formulae $\phi(\vec{x}, y, z)$ and for all $a \in \mathrm{WF}(M)$,
$$\textrm{if } \mathcal{M} \models \exists \vec{x} \phi(\vec{x}, a, b), \textrm{ then } \mathcal{M} \models (\exists \vec{x} \in B) \phi(\vec{x}, a, b).$$
Therefore, by Theorem \ref{Th:MainSelfEmbeddingResult2}, there exists a proper initial self-embedding $j: \mathcal{M} \longrightarrow \mathcal{M}$ such that $j[\mathcal{M}] \subseteq_e B^*$, $j(b)=b$ and $j[\mathcal{M}] \prec_p \mathcal{M}$.
\Square
\end{proof}

\begin{Coroll1} \label{Th:ElementarySelfEmbeddingsForZFCminus}
Let $\mathcal{M}$ be a countable nonstandard model of $\mathrm{ZF}^-+\mathrm{WO}$ such that the well-founded part of $\mathcal{M}$ is contained. Then for all $p \in \omega$ and for all $b \in M$, there exists a proper initial self-embedding $j: \mathcal{M} \longrightarrow \mathcal{M}$ such that $b \in \mathrm{rng}(j)$ and $j[\mathcal{M}] \prec_p \mathcal{M}$.
\end{Coroll1}

Theorem \ref{Th:ExistenceOfElementarySelfEmbeddingsForKPI} also yields the following results that provide two different sufficient conditions for models of $\mathrm{KPI}+\Sigma_{1}\textrm{-Separation}$ to admit proper initial self-embeddings.

\begin{Coroll1} \label{Th:SelfEmbeddingOfExtensionOfKPI1}
Let $\mathcal{M}$ be a countable nonstandard model of $\mathrm{KPI}+\Sigma_{1}\textrm{-Separation}$ such that the well-founded part of $\mathcal{M}$ is contained and let $b \in M$. Then there exists a proper initial self-embedding $j: \mathcal{M} \longrightarrow \mathcal{M}$ such that $b \in \mathrm{rng}(j)$.
\end{Coroll1}

\begin{Coroll1} \label{Th:NonOmegaModelsOfKPISelfEmeddingResult}
Let $\mathcal{M}$ be a countable $\omega$-nonstandard model of $\mathrm{KPI}+\Sigma_{1}\textrm{-Separation}$ and let $b \in M$. Then there exists a proper initial self-embedding $j: \mathcal{M} \longrightarrow \mathcal{M}$ with $b \in \mathrm{rng}(j)$.
\end{Coroll1}

This allows us to give an example of a countable $\omega$-nonstandard model of $\mathrm{MOST}+\Pi_1\textrm{-Collection}$ that admits a proper initial self-embedding, but no proper $\mathcal{P}$-initial self-embedding.

\begin{Examp1} \label{ex:InitialButNoPInitial}
\normalfont Let $\mathcal{M}$ be a countable model of $\mathrm{ZF}+V=L$ that is not $\omega$-standard. Let $\mathcal{N}$ be the substructure of $\mathcal{M}$ with underlying set
$$N= \bigcup_{n \in \omega} (H_{\aleph_n}^{\mathcal{M}})^*.$$
The fact that $\mathcal{N} \models \mathrm{Mac}$ follows immediately from the fact that $\mathcal{N} \models \textrm{Powerset}$ and $\mathcal{N} \subseteq_e^\mathcal{P} \mathcal{M}$. Since $\mathcal{M}$ satisfies the Generalised Continuum Hypothesis, $\mathcal{N} \models \textrm{Axiom }\mathrm{H}$. Therefore, by Lemma \ref{Th:MOSTisMacPlusH}, $\mathcal{N} \models \mathrm{MOST}$. It follows from Lemma \ref{Th:HCutsSatisfyPi1Collection} that $\mathcal{N} \models \mathrm{MOST}+\Pi_1\textrm{-Collection}$. By Corollary \ref{Th:NonOmegaModelsOfKPISelfEmeddingResult}, $\mathcal{N}$ admits a proper initial self-embedding. Now, suppose that $j: \mathcal{N} \longrightarrow \mathcal{N}$ is a proper $\mathcal{P}$-initial self-embedding. But this is impossible, because, since $j[\mathcal{N}] \subseteq_e^\mathcal{P} \mathcal{N}$, cardinals are preserved between $j[\mathcal{N}]$ and $\mathcal{N}$.
\end{Examp1}

We are also able to find an example of a countable $\omega$-nonstandard model of $\mathrm{MOST}+\Pi_1\textrm{-Collection}$ that admits a proper $\mathcal{P}$-initial self-embedding, but no proper rank-initial self-embedding.

\begin{Examp1}
\normalfont Let $\mathcal{N}= \langle N, \in^\mathcal{M} \rangle$ be the countable model of $\mathrm{MOST}+\Pi_1\textrm{-Collection}$ described in Example \ref{ex:InitialButNoPInitial}. Note that $\mathcal{N}$ satisfies the Generalised Continuum Hypothesis and the infinite cardinals of $\mathcal{N}$ are exactly $\aleph_n$ for each standard natural number $n$. In particular, for all $n \in \omega$,
$$\mathcal{N}\models (|\mathcal{P}^n(V_\omega)| = \aleph_n) \textrm{, where } \mathcal{P}^n(V_\omega) \textrm{ is the powerset operation applied to } V_\omega\ n\textrm{-times.}$$
It follows that $\mathcal{N}$ satisfies
\begin{itemize}
\item[] ($\dagger$) For all cardinals $\kappa$, there exists a set $X$ with cardinality $\kappa$ and countable rank.
\end{itemize}
Therefore, $\mathcal{N}$ shows that the theory $\mathrm{MOST}+\Pi_1\textrm{-Collection}+(\dagger)$ is consistent. Now, let $\mathcal{K}= \langle K, \in^\mathcal{K}\rangle$ be a recursively saturated model of $\mathrm{MOST}+\Pi_1\textrm{-Collection}+(\dagger)$. By Theorem \ref{Th:PInitialSelfEmbeddingOfRecursivelySaturatedModels}, $\mathcal{K}$ has a proper $\mathcal{P}$-initial self embedding. Now, suppose $j: \mathcal{K} \longrightarrow \mathcal{K}$ is a proper $\mathcal{P}$-initial self embedding. Since a bijection between an ordinal $\kappa$ and $\alpha \in \kappa$ is a subset of $\mathcal{P}(\kappa \times \kappa)$, for all $\kappa \in \mathrm{rng}(j)$, $\kappa$ is a cardinal according to $\mathcal{K}$ if and only if $\kappa$ is a cardinal according to $j[\mathcal{K}]$. Similarly, if $R \in K$ and $\kappa \in \mathrm{rng}(j)$ is a cardinal of $\mathcal{K}$ such that
$$\mathcal{K} \models (R \subseteq \kappa \times \kappa)\land (R \textrm{ is a well-founded extensional relation with a maximal element}),$$
then $R \in \mathrm{rng}(j)$ and
$$j[\mathcal{K}] \models (R \textrm{ is a well-founded extensional relation with a maximal element})$$
Therefore, by Lemma \ref{Th:ConsequencesOfMOST}, if $\kappa \in \mathrm{rng}(j)$ is a cardinal, then $(\kappa^+)^\mathcal{K} \in \mathrm{rng}(j)$ and $H_{\kappa^+}^{j[\mathcal{K}]}= H_{\kappa^+}^\mathcal{K}$. Now, since $j$ is proper, let $x \in K \backslash \mathrm{rng}(j)$. Let $\kappa \in K$ be such that $\mathcal{K} \models (|\mathrm{TC}(\{x\})|=\kappa)$. By the observations that we have just made, $\kappa \notin \mathrm{rng}(j)$ and for all $y \in K$, if $\mathcal{K} \models (|\mathrm{TC}(\{y\})|\geq \kappa)$, then $y \notin \mathrm{rng}(j)$. Therefore, since $\mathcal{K} \models (\dagger)$, there exists a set $y \in K$ with countable rank in $\mathcal{K}$ such that $y \notin \mathrm{rng}(j)$. This shows that $j$ is not a proper rank-initial self-embedding.
\end{Examp1}

Theorem \ref{Th:ExistenceOfElementarySelfEmbeddingsForKPI} combined with Lemma \ref{Th:WellFoundedPartsOfKPPContained} also yields the following result that shows that every nonstandard model of $\mathrm{KP}^\mathcal{P}+\Sigma_1\textrm{-Separation}$ admits a proper initial self-embedding.

\begin{Coroll1} \label{Th:FriedmanStyleExistenceTheorem}
Let $\mathcal{M}$ be a countable nonstandard model of $\mathrm{KP}^\mathcal{P}+\Sigma_1\textrm{-Separation}$ and let $b \in M$. Then there exists a proper initial self-embedding $j: \mathcal{M} \longrightarrow \mathcal{M}$ with $b \in \mathrm{rng}(j)$.
\end{Coroll1}

Note that the theory $\mathrm{MOST}+\Pi_1\textrm{-Collection}+\Pi_1^\mathcal{P}\textrm{-Foundation}$ is obtained from $\mathrm{KP}^\mathcal{P}+\Sigma_1\textrm{-Separtion}$ by adding the Axiom of Choice. The following example shows that the assumptions of Corollary \ref{Th:FriedmanStyleExistenceTheorem} cannot be weakened to saying that $\mathcal{M}$ is a nonstandard model of $\mathrm{MOST}+\Pi_1\textrm{-Collection}$.

\begin{Examp1}
\normalfont Let $\mathcal{M}$ be a countable model of $\mathrm{ZF}+V=L$ that is $\omega$-standard but has a nonstandard ordinal that is countable according to $\mathcal{M}$. Note that such a model can by obtained from the assumption that there exists a transitive model of $\mathrm{ZF}+V=L$ using \cite[Theorem 2.4]{Keisler-Morley} or, from the same assumption, using the Barwise Compactness Theorem as in the proof of \cite[Theorem 4.5]{mck15}. Let $\mathcal{N}$ be the substructure of $\mathcal{M}$ with underlying set
$$N= \bigcup_{\alpha \in \mathrm{o}(\mathcal{M})} (H_{\aleph_\alpha}^\mathcal{M})^*.$$
Using the same reasoning that was used in Example \ref{ex:InitialButNoPInitial},
$$\mathcal{N} \models \mathrm{MOST}+\Pi_1\textrm{-Collection}.$$
Moreover, $\mathcal{N}$ is nonstandard. Since $\mathcal{M}$ satisfies the Generalised Continuum Hypothesis, a straightforward induction argument inside $\mathcal{M}$ shows that
$$\mathcal{M} \models \forall \alpha(\alpha \textrm{ is an ordinal} \Rightarrow |V_{\omega+\alpha}|= \aleph_\alpha).$$
Therefore,
$$\mathrm{WF}(N)= \bigcup_{\alpha \in \mathrm{o}(\mathcal{N})} (V_{\omega+\alpha}^\mathcal{M})^*$$
and the well-founded part of $\mathcal{N}$ is c-unbounded in $\mathcal{N}$. So, by Theorem \ref{Th:NoSelfEmbeddingWhenStandardPartDense}, $\mathcal{N}$ admits no proper initial self-embedding.
\end{Examp1}

We can also use Theorem \ref{Th:MainSelfEmbeddingResultOmegaNonStandard2} to prove the following variant of Theorem \ref{Th:ExistenceOfElementarySelfEmbeddingsForKPI} for models of extensions of $\mathrm{KPI}$ that are not $\omega$-standard.

\begin{Theorems1}
Let $p \in \omega$, $\mathcal{M}$ be a countable $\omega$-nonstandard model of $\mathrm{KPI}+\Pi_p\textrm{-Collection}+\Pi_{p+2}\textrm{-Foundation}$, and let $b \in M$. Then there exists a proper initial self-embedding $j:\mathcal{M} \longrightarrow \mathcal{M}$ such that $b \in \mathrm{rng}(j)$ and $j[\mathcal{M}] \prec_p \mathcal{M}$.
\end{Theorems1}

\begin{proof}
Consider $\theta(C, n, b, \omega)$ defined by
$$(n \in \omega) \land (\forall m \in n)(\exists \vec{x}\ \mathrm{Sat}_{\Pi_p}(m, \langle \vec{x}, b \rangle) \Rightarrow (\exists \vec{x} \in C) \mathrm{Sat}_{\Pi_p}(m, \langle \vec{x}, b\rangle)).$$
Note that $\Pi_p$-Collection implies that $\theta(C, n, b, \omega)$ is equivalent to a $\Pi_{p+1}$-formula. Moreover, if $n \in \omega$, then there exists a finite set $C$ such that $\theta(C, n, b , \omega)$ holds. Therefore, for all $n \in \omega$,
$$\mathcal{M} \models \exists C \ \theta(C, n, b, \omega).$$
So, by $\Pi_{p+2}$-Foundation, there exists a nonstandard $k \in (\omega^\mathcal{M})^*$ such that
$$\mathcal{M}\models \exists C \ \theta(C, k, b, \omega).$$
Let $C \in M$ be such that
$$\mathcal{M}\models \theta(C, k, b, \omega).$$
And, working inside $\mathcal{M}$, let $B= \mathrm{TC}(C\cup V_\omega)$. Therefore $\mathcal{M} \models \bigcup B \subseteq B$, $\mathrm{WF}(M) \subseteq B^*$ and for all $\Pi_p$-formulae $\phi(\vec{x}, z)$,
$$\textrm{if } \mathcal{M} \models \exists \vec{x} \phi(\vec{x}, b), \textrm{ then } \mathcal{M} \models (\exists \vec{x} \in B) \phi(\vec{x}, b).$$
Therefore, by Theorem \ref{Th:MainSelfEmbeddingResultOmegaNonStandard2}, there exists a proper initial self-embedding $j: \mathcal{M} \longrightarrow \mathcal{M}$ such that $j[\mathcal{M}] \subseteq_e B^*$, $j(b)=b$ and $j[\mathcal{M}] \prec_p \mathcal{M}$. \Square
\end{proof}

\begin{Coroll1} \label{Th:SelfEmbeddingOfExtensionKPI2}
Let $\mathcal{M}$ be a countable $\omega$-nonstandard model of $\mathrm{KPI}+\Pi_2\textrm{-Foundation}$, and let $b \in M$. Then there exists a proper initial self-embedding $j: \mathcal{M} \longrightarrow \mathcal{M}$ such that $b \in \mathrm{rng}(j)$.
\end{Coroll1}

Note that Corollaries \ref {Th:NonOmegaModelsOfKPISelfEmeddingResult} and \ref{Th:SelfEmbeddingOfExtensionKPI2} give two distinct extensions of $\mathrm{KPI}$ such that every countable $\omega$-nonstandard model of these extensions is isomorphic to a transitive proper initial segment of itself.

We will next use Theorem \ref{Th:MainSelfEmbeddingResultOmegaNonStandard2} together with Corollary \ref{cor. of Levy-Shoenfield} to verify the surprising result that every model of $\mathrm{ZFC}$ that is $\omega$-nonstandard is isomorphic to a transitive substructure of the hereditarily countable sets of its own $L$.

\begin{Theorems1} \label{Th:InitialEmbeddingIntoL}
Let $\mathcal{M}$ be a countable $\omega$-nonstandard model of $\mathrm{ZF}$. Then there exists a proper initial self-embedding $j: \mathcal{M} \longrightarrow \mathcal{M}$ such that $j[\mathcal{M}] \subseteq_e (H_{\kappa}^{L^\mathcal{M}})^*$, where $\kappa= (\aleph_1^L)^\mathcal{M}$.
\end{Theorems1}

\begin{proof}
Let $\kappa= (\aleph_1^L)^\mathcal{M}$. Now, let $B= H_\kappa^{L^\mathcal{M}}$. It is clear that $\mathcal{M} \models \bigcup B \subseteq B$ and $\mathrm{WF}(M) \subseteq B^*$. Note that $\emptyset \in B^* \cap M$. By Corollary \ref{cor. of Levy-Shoenfield}, for all $\Delta_0$-formula $\phi(\vec{x}, z)$,
$$\textrm{if } \mathcal{M} \models \exists \vec{x} \phi(\vec{x}, \emptyset), \textrm{ then } \mathcal{M} \models (\exists \vec{x} \in B) \phi(\vec{x}, \emptyset).$$
So, by Theorem \ref{Th:MainSelfEmbeddingResultOmegaNonStandard2}, there exists a proper initial self-embedding $j: \mathcal{M} \longrightarrow \mathcal{M}$ such that $j[\mathcal{M}] \subseteq_e B^*$.
\Square
\end{proof}

Hamkins \cite{ham13} showed that if $\mathcal{M}$ is a countable model of $\mathrm{ZF}$, then there exists an embedding of $\mathcal{M}$ into its own $L$. However, the embeddings produced in \cite{ham13} are not required to be initial embeddings. Theorem \ref{Th:InitialEmbeddingIntoL} shows that under the condition that $\mathcal{M}$ is a countable $\omega$-nonstandard model of $\mathrm{ZF}$, there exists an embedding of $\mathcal{M}$ into its own $L$ that is also initial.  Question 35 of \cite{ham13} asks whether every countable model of set theory can be embedded into its own $L$ by an embedding that preserves ordinals. Since initial embeddings preserve ordinals, Theorem \ref{Th:InitialEmbeddingIntoL} provides a positive answer to this question when $\mathcal{M}$ is a countable $\omega$-nonstandard model of $\mathrm{ZF}$.

Theorem \ref{Th:InitialEmbeddingIntoL} immediately implies the corollary below that shows that every countable model of $\mathrm{ZF}$ that is not $\omega$-standard can be end-extended to a model of $\mathrm{ZFC}+V=L$.

\begin{Coroll1} \label{Th:SpecialCaseOfBarwiseResult}
Let $\mathcal{M}$ be a countable $\omega$-nonstandard model of $\mathrm{ZF}$. Then there exists structures $\mathcal{N}_1$ and $\mathcal{N}_2$ such that
\begin{itemize}
\item[(I)] $\mathcal{M} \subseteq_e \mathcal{N}_1 \subseteq_e \mathcal{N}_2$,
\item[(II)] $\mathcal{N}_2 \models \mathrm{ZFC}+V=L$, and
\item[(III)] $\mathcal{N}_1= \langle (H_{\aleph_1}^{\mathcal{N}_2})^*, \in^{\mathcal{N}_2} \rangle$.
\end{itemize}
\end{Coroll1}


Corollary \ref{Th:SpecialCaseOfBarwiseResult} is a special case of \cite[Theorem 3.1]{bar71}, which shows that Corollary \ref{Th:SpecialCaseOfBarwiseResult} holds for all countable models of $\mathrm{ZF}$. Barwise used methods from infinitary logic; Hamkins has recently formulated a purely set-theoretic proof of the same result \cite{ham18}.

We now turn to applying Theorem \ref{Th:MainSelfEmbeddingResult2} to finding transitive partially elementary substructures of nonstandard models of $\mathrm{ZF}^-+\mathrm{WO}$. Despite the failure of reflection in $\mathrm{ZF}^-+\mathrm{WO}$ \cite{fgk19}, Quinsey \cite[Corollary 6.9]{qui80} employed indicators and methods from infinitary logic to show the following:

\begin{Theorems1} \label{Th:QuinseyResult}
(Quinsey) Let $n \in \omega$. If $\mathcal{M} \models \mathrm{ZF}^-$, then there exists $\mathcal{N} \subseteq_e \mathcal{M}$ such that $N \neq M$, $\mathcal{N} \prec_n \mathcal{M}$ and $\mathcal{N} \models \mathrm{ZF}^-$.
\end{Theorems1}

Extensions of H. Friedman's self-embedding result \cite{fri73} proved by Gorbow \cite{gor18} show that if the nonstandard model $\mathcal{M}$ in Theorem \ref{Th:QuinseyResult} is countable and satisfies $\mathrm{ZF}$, then the conclusions Theorem \ref{Th:QuinseyResult} can be strengthened to require that $\mathcal{N} \subseteq_e^\mathcal{P} \mathcal{M}$ and $\mathcal{N} \cong \mathcal{M}$. In light of this, it natural to ask under what circumstances the conclusion of Theorem \ref{Th:QuinseyResult} can be strengthened to require that the $\Sigma_n$-elementary submodel be isomorphic to the original nonstandard model. Theorem \ref{Th:ModelOfZFCminusWithoutSelfEmbedding} shows that such a strengthening of Theorem \ref{Th:QuinseyResult} does not hold in general, even when $n=0$ and the model is countable. However, using Theorem \ref{Th:MainSelfEmbeddingResult2}, we can show that the countable nonstandard models of $\mathrm{ZF}^-+\mathrm{WO}+\forall \alpha\ \Pi_\infty^1-\mathrm{DC}_\alpha$ for which this strengthening of Quinsey's result holds are exactly the models in which the well-founded part is c-bounded.

Our final result below shows that the c-boundedness of the well-founded part of a countable nonstandard model of $\mathcal{M}$ of $\mathrm{ZF}^-+\mathrm{WO}+\forall \alpha\ \Pi^1_\infty\mathrm{-DC}_\alpha$ is necessary and sufficient for $\mathcal{M}$ to admit a proper initial self-embedding.

\begin{Theorems1} \label{Th:IsomorphicElementarySubmodelsResult}
Let $\mathcal{M}$ be a countable nonstandard model of $\mathrm{ZF}^-+\mathrm{WO}+\forall \alpha\ \Pi^1_\infty\mathrm{-DC}_\alpha$. Then the following are equivalent:
\begin{itemize}
\item[(I)] The well-founded part of $\mathcal{M}$ is c-bounded in $\mathcal{M}$,
\item[(II)] $\mathrm{WF}(\mathcal{M}) \models \mathrm{Powerset}$,
\item[(III)] For all $n \in \omega$ and for all $b \in M$, there exists a proper initial self-embedding $j: \mathcal{M} \longrightarrow \mathcal{M}$ such that $b \in \mathrm{rng}(j)$ and $j[\mathcal{M}] \prec_n \mathcal{M}$.
\end{itemize}
\end{Theorems1}

\begin{proof}
(I)$\Rightarrow$(II) is Lemma \ref{Th:CBoundedImpliesPowersetInWellFoundedPart}. To see that (II)$\Rightarrow$(III), assume that (II) holds and note that Lemma \ref{Th:PowersetInWellFoundedPartImpliesContained} implies that the standard part of $\mathcal{M}$ is contained. Therefore, by Corollary \ref{Th:ElementarySelfEmbeddingsForZFCminus}, (III) holds. Finally, (III)$\Rightarrow$(I) is the contrapositive of Theorem \ref{Th:NoSelfEmbeddingWhenStandardPartDense}. \Square
\end{proof}

\section[Questions]{Questions} \label{Sec:questions}

\noindent \textbf{Question 6.1} \textit{Can Theorem} \ref{Th:ExistenceOfElementarySelfEmbeddingsForKPI} \textit{be strengthened by adding the requirement to the conclusion of the theorem that the self-embedding} $j$ \textit{fixes every member of} $b^*$?
\begin{itemize}
\item The above question is motivated by the ``moreover" clause of Theorem \ref{Th_Gorbow}.
\end{itemize}

\noindent \textbf{Question 6.2} \textit{Does every countable model of} $\mathrm{KPI}$ \textit{that is not} $\omega$-\textit{standard admit a proper initial self-embedding}?

\begin{itemize}
\item The above question is prompted by Corollaries \ref{Th:SelfEmbeddingOfExtensionOfKPI1} and \ref{Th:SelfEmbeddingOfExtensionKPI2}, which provide sufficient conditions for a countable model of $\mathrm{KPI}$ that is not $\omega$-standard to admit a proper initial self-embedding.
\end{itemize}

\noindent \textbf{Question 6.3} \textit{Is there a countable model of} $\mathcal{M}$ \textit{of} $\mathrm{ZF}^-+\mathrm{WO}$ \textit{such that the well-founded part of} $\mathcal{M}$ i\textit{s c-bounded} \textit{and} $\mathcal{M}$ \textit{does not admit any proper initial self-embedding?}
\begin{itemize}
\item The key role played by the scheme $\forall \alpha\ \Pi^1_\infty\mathrm{-DC}_\alpha$ in the proof of Theorem \ref{Th:IsomorphicElementarySubmodelsResult} suggests a positive answer to the above question.

\end{itemize}

\bibliographystyle{alpha}
\bibliography{.}

\begin{thebibliography}{9}


\bibitem[BE]{BahEna18} Bahrami, Saeideh; and Enayat, Ali. ``Fixed points of self-embbeddings of models of arithmetic". \emph{Annals of Pure and Applied Logic}. Vol. 168. 2018. pp 487--513.


\bibitem[Bar71]{bar71} Barwise, Jon. ``Infinitary methods in the model theory of set theory". \emph{Logic Colloquium '69 (Proceedings Summer School and Colloquium, Manchester, 1969)}. North-Holland, Amsterdam. 1971. pp 53--66.

\bibitem[Bar75]{bar75} Barwise, Jon. \emph{Admissible Sets and Structures}. Perspectives in Mathematical Logic. Springer-Verlag. 1975.

\bibitem[BF]{barwise-fischer70} Barwise, Jon; and Fischer, Edward. ``The Shoenfield absolutness lemma". \emph{Israel J. Math}. Vol. 8. 1970. pp. 329--339.



\bibitem[CK]{Chang-Keisler}Chang, Chen C.; and Keisler, H. Jerome. \emph{Model Theory}, Studies in Logic and the Foundations of Mathematics, 3rd ed., North-Holland publishing Co., Amsterdam. 1990.


\bibitem[EKM]{ekm18} Enayat, Ali; Kaufmann, Matt; and McKenzie, Zachiri. ``Largest initial segments pointwise fixed by automorphisms of models of set theory". \emph{Archive for Mathematical Logic}. Vol. 57. No. 1-2. 2018. pp 91--139.

\bibitem[Fla]{fla75} Flannagan, Timothy B. ``Axioms of choice in Morse-Kelley class theory". {\it $\vDash$ISILC Logic Conference}. Edited by G. H. M\"{u}ller, A. Oberschelp, K. Potthoff. Springer Lecture Notes in Mathematics. Vol. 499. Springer, Berlin, Heidelberg. 1975. pp 190--247.

\bibitem[FK]{fk91} Forster, Thomas; and Kaye, Richard. ``End-extensions preserving power set". \emph{The Journal of Symbolic Logic}. Vol. 56. No. 1. 1991. pp 323--328.

\bibitem[Fri]{fri73} Friedman, Harvey M. ``Countable models of set theories". \emph{Cambridge Summer School in Mathematical Logic, August 1--21, 1971}. Edited by A. R. D. Mathias and H. Rogers Jr. Springer Lecture Notes in Mathematics. Vol. 337. Springer, Berlin. 1973. pp 539--573.

\bibitem[FLW]{flw16} Friedman, Sy-David; Li, Wei; and Wong, Tin Lok. ``Fragments of Kripke-Platek Set Theory and the Metamathematics of $\alpha$-Recursion Theory". \emph{Archive for Mathematical Logic}. Vol. 55. No. 7. 2016. pp 899--924.

\bibitem[FGK]{fgk19} Friedman, Sy-David; Gitman, Victoria; and Kanovei, Vladimir. ``A model of second-order arithmetic satisfying AC but not DC". {To appear in the Journal of Mathematical Logic}.

\bibitem[GHT]{ght16} Gitman, Victoria; Hamkins, Joel D.; and Johnstone, Thomas A. ``What is the theory $\mathrm{ZFC}$ without powerset". \emph{Mathematical Logic Quarterly}. Vol. 62. No. 4-5. 2016. pp 391--406.

\bibitem[Gor]{gor18} Gorbow, Paul Kindvall. \emph{Self-similarity in the Foundations}, Doctoral Dissertation, University of Gothenburg, 2018.
Available online: {\tt https://arxiv.org/abs/1806.11310}.

\bibitem[Ham13]{ham13} Hamkins, Joel D. ``Every countable model of set theory embeds into its own constructible universe". \emph{Journal of Mathematical Logic}. Vol. 13. 2013.

\bibitem[Ham18]{ham18} Hamkins, Joel D. ``A new proof of the Barwise extension theorem, without infinitary logic". CUNY Logic Workshop, December 2018, Available online:
{\tt http://jdh.hamkins.org}.

\bibitem[J]{Jechbook} Jech, Thomas. \emph{Set Theory}. Springer Monographs in Mathematics, Springer, Berlin (2003).

\bibitem[KW] {kw07} Kaye, Richard; and Wong Tin Lok. ``On interpretations of arithmetic and set theory''. \emph{Notre Dame Journal of Formal Logic}. Vol. 48. 2007. pp. 497--510.

\bibitem[KM]{Keisler-Morley} Keisler, H. Jerome; and Morley, Michael. ``Elementary
extensions of models of set theory". \emph{Israel J. Math}. Vol. 5. 1968.  pp 49--65.


\bibitem[L\'{e}v]{lev64} L\'{e}vy, Azriel. ``The interdependence of certain consequences of the axiom of choice". \emph{Fundamenta Mathematicae}. Vol. 54. 1964. pp 135--157.

\bibitem[McK15]{mck15} McKenzie, Zachiri. ``Automorphisms of models of set theory and extensions of $\mathrm{NFU}$". \emph{Annals of Pure and Applied Logic}. Vol. 166. 2015. pp 601--638.

\bibitem[McK19]{mck19} McKenzie, Zachiri. ``On the relative strengths of fragments of collection". \emph{Mathematical Logic Quarterly}. Vol. 65. No. 1. 2019. pp 80--94.

\bibitem[Mat69]{mat69} Mathias, Adrian R. D. ``Notes on set theory". Available online: {\tt https://www.dpmms.ac.uk/\textasciitilde ardm/}.

\bibitem[Mat01]{mat01} Mathias, Adrian R. D. ``The strength of Mac Lane set theory". \emph{Annals of Pure and Applied Logic}. Vol. 110. 2001. pp 107--234.

\bibitem[Qui]{qui80} Quinsey, Joseph E. ``Some problems in logic". {\it Ph.D. Thesis}. University of Oxford. 1980.

\bibitem[Res]{res87} Ressayre, Jean-Pierre. ``Mod\`{e}les non standard et sous-syst\`{e}mes remarquables de ZF". \emph{Mod\`{e}les non standard en arithm\'{e}tique et th\'{e}orie des ensembles}. Publications Math\'{e}matiques de l'Universit\'{e} Paris VII. Vol. 22. Universit\'{e} de Paris VII, U.E.R. de Math\'{e}matiques, Paris. 1987. pp 47--147.



\bibitem[Tak]{tak72} Takahashi, Moto-o. ``$\tilde{\Delta}_1$-definability in set theory". \emph{Conference in mathematical logic --- London '70}. Edited by W. Hodges. Springer Lecture Notes in Mathematics. Vol. 255. Springer. 1972. pp 281--304.

\bibitem[Zar82]{zar82} Zarach, Andrzej M. ``Unions of $\mathrm{ZF}^-$-models which are themselves $\mathrm{ZF}^-$-models". {\it Logic Colloquium '80 (Prague, 1980), Studies in  Logic and the Foundations of Mathematics}. Vol. 108. North-Holland, Amsterdam. 1982. pp
315--342.

\bibitem[Zar96]{zar96} Zarach, Andrzej M. ``Replacement $\nrightarrow$ Collection". \emph{G\"{o}del '96 (Brno, 1996)}. Vol. 6. Berlin: Springer. 1996. pp 307--322.

\end{thebibliography}

\end{document}